\documentclass[12pt]{article}

\usepackage{amssymb}   
\usepackage{amsthm}    
\usepackage{amsmath}   
\usepackage{stmaryrd}  
\usepackage{titletoc}  
\usepackage{mathrsfs}  

\newlength{\defbaselineskip}
\newcommand{\setlinespacing}[1]%
           {\setlength{\baselineskip}{#1 \defbaselineskip}}

\theoremstyle{plain}
\newtheorem{thm}{Theorem}[section]
\newtheorem{cor}[thm]{Corollary}
\newtheorem{lem}[thm]{Lemma}

\theoremstyle{definition}
\newtheorem{defn}{Definition}[section]
\newtheorem{ass}{Assumption}
\newtheorem{rmk}{Remark}[section]

\newcommand{\li}{\left|\!\left|\!\left|}
\newcommand{\ri}{\right|\!\right|\!\right|}
\newcommand{\la}{\langle}
\newcommand{\ra}{\rangle}

\makeatletter
   
   \@addtoreset{equation}{section}
\makeatother

\begin{document}

\title{On the Dirichlet Problem for Backward Parabolic Stochastic
Partial Differential Equations\\ in General Smooth Domains\footnotemark[1]}

\author{Kai Du\footnotemark[2] \and Shanjian Tang\footnotemark[2]}
\date{October 10, 2009}

\footnotetext[1]{Supported by NSFC Grant \#10325101, Basic Research Program
of China (973 Program)  Grant \# 2007CB814904.}

\footnotetext[2]{Department of Finance and Control Sciences, School of
Mathematical Sciences, Fudan University, Shanghai 200433, China. E-mail:
kdu@fudan.edu.cn (Kai Du), sjtang@fudan.edu.cn (Shanjian Tang).}

\maketitle


\abstract{Backward stochastic partial differential equations of parabolic
type with variable coefficients are considered in smooth domains. Existence
and uniqueness results are given in weighted Sobolev spaces allowing the
derivatives of the solutions to blow up near the boundary.}
\medskip

\medskip

{\bf Abbreviated Title.} On the Dirichlet problem for BSPDEs
\medskip



\section{Introduction}

In this paper we consider the Dirichlet problem of backward stochastic
partial different equations (BSPDEs, for short) of the form
\begin{eqnarray}\label{eq:a1}
\begin{split}
d p (t,x)=&-\bigg{[}\frac{\partial}{\partial x^i}\bigg{(}
a^{ij}(t,x)\frac{\partial p}{\partial x^j}(t,x) +\sigma^{ik}(t,x)
q^{k}(t,x)\bigg{)} +b^i(t,x)\frac{\partial p}{\partial
x^i}(t,x)-c(t,x)p(t,x)\\
&+\nu^k(t,x)q^k(t,x)+F(t,x)\bigg{]}dt+q^k(t,x)dW^k_t,\quad (t,x)\in
[0,T]\times \cal{D},
\end{split}\end{eqnarray}
where $\mathcal{D}$ is a domain of $d$-dimensional Euclidean space, where $W
\triangleq \{W^{k}_{t};t\geq 0\}$ is a $d_1$-dimensional Wiener process
generating a natural filtration $\{\mathscr{F}_{t}\}_{t\geq 0}$, and the
functions $p(t,x)=p(\omega,t,x)$ and $q(t,x)=q(\omega,t,x)$ are both unknown
of which the first one should satisfy given terminal conditions on
$\Omega\times \{T\}\times\mathcal{D}$ and boundary conditions on
$\Omega\times (0,T)\times\partial{\mathcal{D}}$. The coefficients
$a,b,c,\sigma,\nu$ and the free term $F$ are random functions depending on
$(t,x)$. An adapted solution of this equation is a $\mathscr{P}\times
B(\mathcal{D})$-measurable function pair $(p,q)$ satisfying the given
terminal-boundary conditions and solving Eq.\eqref{eq:a1} under some
appropriate sense, where $\mathscr{P}$ is the $\sigma$-algebra of predictable
sets on $\Omega \times (0,T)$ associated with $\{\mathscr{F}_{t}\}_{t\geq
0}$.

BSPDEs, which is a natural extension of backward SDEs (see e.g. Pardoux-Peng
\cite{PaPe90} or El Karoui-Peng-Quenez \cite{KPQ97}), arise in many
applications of probability theory and stochastic processes, for instance in
the optimal control of processes with incomplete information, as adjoint
equations of Duncan-Mortensen-Zakai filtration equations (see e.g. Bensoussan
\cite{Bens92}, Nagasa-Nisio \cite{NaNi90}, Zhou \cite{Zhou93} and Tang
\cite{Tang98a,Tang98b}). A class of fully nonlinear BSPDEs, the so-called
backward stochastic Hamilton-Jacobi-Bellman equations, are also introduced in
the study of controlled non-Markovian processes by Peng \cite{Peng92a}. For
more aspects of BSPDEs, we refer to e.g. Barbu-R\u{a}\c{s}canu-Tessitore
\cite{BRT03}, Ma-Yong \cite{MaYo97} and Tang-Zhang \cite{TaZh09}.

The works concerning the existence and uniqueness of the adapted solution to
a BSPDE can be found in e.g. Bensoussan \cite{Bens92}, Hu-Peng \cite{HuPe91},
Peng \cite{Peng92a}, Zhou \cite{Zhou92}, Tessitore \cite{Tess96}, Ma-Yong
\cite{MaYo99} and Tang \cite{Tang05}. However, most of these results are
concerned with the Cauchy problem for BSPDEs, i.e. the case of
$\mathcal{D}=\mathbb{R}^{d}$. Indeed, the existence and uniqueness of the
weak solution (see Definition \ref{defn:b1}) of Eq.\eqref{eq:a1} in a general
domain can be easily obtained using the finite-dimensional approximation
method (see e.g. \cite{Zhou92}). However, the weak solution usually has poor
smoothness properties. It is interesting to know the regularity with respect
to $x$ of the weak solution in association with the coefficients and free
term of Eq.\eqref{eq:a1} and the boundary data. The issue was investigated by
Zhou \cite{Zhou92}, Ma-Yong\cite{MaYo99}, etc. for the Cauchy problem. Their
approaches strictly depend on the fact that the derivatives of the unknown
functions $p,q$ with respect to $x$ vanish at infinity, which do not work for
the Dirichlet problem since the derivatives of $p$ on the boundary is unknown
and seems to be very difficult to handle. Note that a special form of
Eq.\eqref{eq:a1} with Dirichlet conditions was studied in \cite{Tess96} by
using the method of semigroups, which required that the coefficients were
independent of $(\omega,t)$. To our best knowledge, there was no general
result which concerns the Dirichlet problem for BSPDE \eqref{eq:a1}.

In this paper, we follow the approach of Krylov \cite{Kryl94} and establish
the existence and uniqueness results for the Dirichlet problem of
Eq.\eqref{eq:a1} in weighted Sobolev spaces allowing the derivatives of the
solution to blow up near the boundary of $\mathcal{D}$ and concluding the
interior regularity of the solution. The requirements on the coefficients are
rather weak. The main difference from the investigation of SPDEs in
\cite{Kryl94} is that the (adapted) solution of \eqref{eq:a1}, containing two
components, has no explicit formulation without involving conditional
expectation even in the simplest situation, which prevents us to directly
prove the regularity of solutions of 1-dimensional equations. We shall
overcome this difficulty with the aid of the method of continuation and the
difference quotient. Just as the existing results for the Cauchy problem of
BSPDEs, our work is also in the framework of $W^n_2$. Indeed, it seems to be
a big challenge to establish an $L^p$-theory ($p\geq 2$) for BSPDEs. For the
$L^p$-theory of SPDEs, we refer to e.g. Krylov \cite{Kryl99a}, Kim-Krylov
\cite{KiKr04a} and Kim \cite{Kim04}.

This paper is organized as follows. In Section 2, we give some preliminaries
and present our main results. Section 3 and Section 4 are the most essential
parts of the whole paper, in which we investigate the equation on the half
axis. In Section 5, we obtain a result for equations in the half space with
the coefficients independent of $x$. This result help us, in Section 6, to
proof a theorem which deals with equations in the half space. In Section 7,
we complete the proof of our main theorem with the help of the result in
Section 6.

\section{Main results}

Let $(\Omega,\mathscr{F},\{\mathscr{F}_{t}\}_{t\geq 0},P)$ be a complete
filtered probability space on which is defined a $d_1$-dimensional Wiener
process $W=\{W_t;t\geq 0\}$ such that $\{\mathscr{F}_{t}\}_{t\geq 0}$ is the
natural filtration generated by $W$, augmented by all the $P$-null sets in
$\mathscr{F}$. Fix a positive number $T$. Denote by $\mathscr{P}$ the
$\sigma$-algebra of predictable sets on $\Omega \times (0,T)$ associated with
$\{\mathscr{F}_{t}\}_{t\geq 0}$.

Let $\mathcal{D}$ be a domain in $\mathbb{R}^{d}$ with boundary of class
$C^{n+2}$, where $n$ is a nonnegative integer.

For the sake of convenience, we denote $$D_{i}=\frac{\partial}{\partial
x^i},\quad D_{ij}=\frac{\partial^2}{\partial x^i \partial x^j},\quad
i,j=1,\dots,d,$$ and for any multi-index $\alpha=(\alpha_1,\dots,\alpha_d)$
$$D^{\alpha}=\bigg{(} \frac{\partial}{\partial x^1} \bigg{)}^{\alpha_1}
\bigg{(} \frac{\partial}{\partial x^2} \bigg{)}^{\alpha_2} \cdots \bigg{(}
\frac{\partial}{\partial x^d} \bigg{)}^{\alpha_d}, \quad
|\alpha|=\alpha_1+\cdots +\alpha_d.$$ Moreover, denote by $Du$ and $D^{2}u$
respectively the gradient and the Hessian matrix for the function $u$ defined
on $\mathbb{R}^{d}$. We will also use the summation convention.

Throughout the paper, by saying that a vector-valued or matrix-valued
function belongs to a function space (for instance, $Du\in
L^{2}(\mathcal{D})$), we mean all the components belong to that space.

Denote
\begin{eqnarray*}\begin{split}
& \mathbb{H}^{0}(\mathcal{D})=L^{2}(\Omega\times (0,T),
\mathscr{P},L^{2}(\mathcal{D})),\\
& \mathbb{H}^{1}_{0}(\mathcal{D})=L^{2}(\Omega\times (0,T),
\mathscr{P},H^{1}_{0}(\mathcal{D})),\\
& \mathbb{H}^{-1}(\mathcal{D})=L^{2}(\Omega\times (0,T),
\mathscr{P},H^{-1}(\mathcal{D})). \end{split}
\end{eqnarray*}
\smallskip

Let $V$ and $H$ be two separable Hilbert spaces such that $V$ is densely
edmbedded in $H$. We identify $H$ with its dual space, and denote by $V'$ the
dual of $V$. We have then $V \subset H \subset V'$. Denote by $\|\cdot\|_{H}$
the norms of $H$, by $\la\cdot,\cdot\ra_{H}$ the scalar product in $H$, and
by $\la\cdot,\cdot\ra$ the duality product between $V$ and $V'$.

We consider three processes $v(\cdot,\cdot),m(\cdot,\cdot)$ and
$v'(\cdot,\cdot)$ defined on $\Omega\times[0,T]$ with values in $V,H$ and
$V'$ respectively. Let $v(\omega,t)$ be measurable with respect to
$(\omega,t)$ and be $\mathscr{F}_{t}$-measurable with respect to $\omega$ for
a.e. $t\in[0,T]$; for any $\eta\in V$ the quantity $\la\eta, v'(\omega,t)\ra$
is $\mathscr{F}_{t}$-measurable in $\omega$ for a.e. $t$ and is measurable
with respect to $(\omega,t)$. Assume that $m(\omega,t)$ is strongly
continuous in $t$ and is $\mathscr{F}_{t}$-measurable with respect to
$\omega$ for any $t$, and is a local martingale. Let $\la m \ra$ be the
increasing process for $\|m\|_{H}^{2}$ in the Doob-Meyer Decomposition (see
e.g. \cite[Page 1240]{KrRo81}).

Proceeding identically to the proof of Theorem 3.2 in Krylov-Rozovskii
\cite{KrRo81}, we have the following result concerning It\^o's formula, which
is the backward version of \cite[Thm 3.2]{KrRo81}.

\begin{lem}\label{lem:b5}
  Let $\varphi\in L^{2}(\Omega,\mathscr{F}_{T},H)$. Suppose that for any
  $\eta \in V$ and any $t\in [0,T]$, it holds almost
  surely that
  \begin{equation*}
    \la \eta,v(t)\ra_{H} = \la \eta,\varphi \ra_{H}
    + \int_{t}^{T} \la \eta,v'(s) \ra ds
    + \la \eta, m(T)-m(t) \ra_{H}.
  \end{equation*}
  Then there exist a set $\Omega'\subset \Omega$ s.t. $P(\Omega')=1$
  and a function $h(t)$ with values in $H$ such that
  \begin{enumerate}
    \item[a)] $h(t)$ is $\mathscr{F}_{t}$-measurable for any $t\in [0,T]$ and
    strongly continuous with respect to
    $t$ for any $\omega$, and $h(t)=v(t)$ (in the space $H$) a.s.
    $(\omega,t)\in \Omega\times [0,T]$, and $h(T)=\varphi$ for any
    $\omega\in\Omega'$;
    \item[b)] for any $\omega\in \Omega'$ and any $t\in [0,T]$,
    \begin{equation*}
      \|h(t)\|_{H}^{2} = \|\varphi\|_{H}^{2} + 2\int_{t}^{T}\la v(s),v'(s)\ra
      ds + 2 \int_{t}^{T}\la h(s), d m(s)\ra_{H} + \la m \ra_{T} - \la m
      \ra_{t}.
    \end{equation*}
  \end{enumerate}
\end{lem}
\medskip

Now consider Eq.\eqref{eq:a1} with the Dirichlet boundary condition
\begin{equation}\label{con:b}
\left\{\begin{array}{l}
p(t,x)=0,\ t\in [0,T],\ x\in\partial \mathcal{D},\\
p(T,x)=\phi(x),\ x\in \mathcal{D}.
\end{array}\right.
\end{equation}

\begin{defn}\label{defn:b1}
A $\mathscr{P}\times B(\mathcal{D})$-measurable function pair $(p,q)$ valued
in $\mathbb{R}\times \mathbb{R}^{d_1}$ is called a (generalized or weak)
solution of Eq.\eqref{eq:a1}, if $p\in \mathbb{H}^{1}_{0}(\mathcal{D})$ and
$q \in \mathbb{H}^{0}(\mathcal{D})$, such that for any $\eta \in
H_{0}^{1}(\mathcal{D})$ and any $t\in [0,T]$, it holds almost surely that
\begin{eqnarray*}
\begin{split}
\int_{\mathcal{D}}p(t,x)\eta(x)dx = & \int_{\mathcal{D}}\phi(x)\eta(x)dx
+\int_{t}^{T}\int_{\mathcal{D}}\bigg{\{}
D_{i}\big[a^{ij}(t,x)D_{j}p(t,x)+\sigma^{ik}(t,x)q^{k}(t,x)\big]\\
&+b^{i}(t,x)D_{i}p(t,x) -c(t,x)p(t,x)+\nu^k(t,x)q^k(t,x)\\
&+F(t,x)\bigg{\}}\eta(x)dx dt
-\int_{t}^{T}\int_{\mathcal{D}}q^k(t,x)\eta(x) dx dW^k_t.
\end{split}
\end{eqnarray*}
\end{defn}

Such a definition can be found in \cite{MaYo99} for the Cauchy problem of
BSPDEs. The weak solution for SPDEs is referred to \cite{KrRo81}.

\begin{ass}\label{ass:b1}
The given functions $a,b,c,\sigma,\nu$ and $F$ are $\mathscr{P} \times
B(\mathcal{D})$-measurable with values in the set of real symmetric $d\times
d$ matrices, in $\mathbb{R}^{d}$, in $\mathbb{R}$, in $\mathbb{R}^{d\times
d_{1}}$, in $\mathbb{R}^{d_1}$ and in $\mathbb{R}$ respectively. The real
function $\phi$ is $\mathscr{F}_T\times B(\mathcal{D})$-measurable.
\end{ass}

Using the duality method as in Zhou \cite{Zhou92}, in view of Lemma
\ref{lem:b5}, we can prove the following

\begin{lem}\label{lem:b3}
Let the functions $a,b,c,\sigma,\nu$ satisfy Assumption \ref{ass:b1} and be
bounded by $K$. Assume that $\kappa I+ \sigma\sigma^{*}\leq 2a \leq
\kappa^{-1}I$. Suppose $F\in \mathbb{H}^{-1}(\mathcal{D}),\; \phi \in
L^{2}(\Omega,\mathscr{F}_{T},L^{2}(\mathcal{D}))$. Then Eq.\eqref{eq:a1} with
the boundary condition \eqref{con:b} has a unique solution $(p,q)$ in the
space $\mathbb{H}^{1}_{0}(\mathcal{D})\otimes\mathbb{H}^{0}(\mathcal{D})$
such that
\begin{eqnarray*}
E\sup_{t\leq T}\|p(t,\cdot)\|_{L^{2}(\mathcal{D})}^{2} + \| D p \|
_{\mathbb{H}^{0}(\mathcal{D})}^{2} +\| q\|_{\mathbb{H}^{0}(\mathcal{D})}^{2}
\leq C\big{(} \| F \|_{\mathbb{H}^{-1}(\mathcal{D})}^{2} + \|
\phi\|_{L^{2}(\Omega,L^{2}(\mathcal{D}))}^{2}\big{)},
\end{eqnarray*}
where the constant $C=C(\kappa,K,T)$.
\end{lem}
A quite related result can be found in Peng \cite[Thm 2.2]{Peng92a}. The
counterpart result for SPDEs is referred to \cite{KrRo81} and \cite{Rozo90}.
\bigskip

In this paper, we investigate the regularity of the (weak) solution of
Eq.\eqref{eq:a1} under some appropriate conditions on the coefficients, the
free term and the boundary data. Before stating our main result, we introduce
some more notations and assumptions.

Fix some constants $K\in (1,\infty)$ and $\rho_{0},\kappa\in (0,1)$.

As well as the smoothness of $\partial{\mathcal{D}}$, we assume that there is
a function $\psi \in C^{n+2}(\bar{\mathcal{D}})$ such that $\psi(x)\leq
\kappa^{-1} \textrm{dist}(x,\partial \mathcal{D})$ if
$\textrm{dist}(x,\partial \mathcal{D})\leq \rho_{0}$ and $\psi(x)\geq \kappa
(\rho_{0} \wedge \textrm{dist}(x,\partial \mathcal{D}))$ on $\mathcal{D}$.
For instance, the function $\psi$ exists for any bounded domain as a result
of the Heine-Borel Theorem.

Denote $B_{+}=\{x\in \mathbb{R}^{d}:|x|<1,x^1>0\}, B_{\rho}(x)=\{y\in
\mathbb{R}^{d}:|x-y|<\rho \}$.

\begin{ass}\label{ass:b2}
For every $x\in \partial \mathcal{D}$ there exist a domain $U\subset
B_{8K\rho_0}(x)$ containing the ball $B_{4\rho_0}(x)$ and a one-to-one map
$\Phi:2B_{+}\rightarrow U\cap \mathcal{D}$ having the properties:
$$x=\Phi(0),\quad\Phi(B_{+})\supset B_{4\rho_0}(x) \cap \mathcal{D},\quad
\kappa |\xi|^2\leq \big{|}(D \Phi) \xi\big{|}^{2}\leq \kappa^{-1}|\xi|^2$$
for any $\xi\in \mathbb{R}^{d}$, and $|D^{\alpha}\Phi|\leq K$ for any
multi-index $\alpha$ with $|\alpha|\leq n+2$. Here $D\Phi$ is the Jacobi
matrix of $\Phi$.

As well as $\psi$ the map $\Phi$ exists for any bounded domain (with an
appropriate $\rho_{0}$).
\end{ass}

\begin{ass}\label{ass:b3}
For any multi-indices $\alpha,\beta$ such that $|\alpha|\leq n$ and
$|\beta|\leq n+1$, we have
\begin{eqnarray*}
\begin{split}
&\kappa I+\sigma\sigma^{*}\leq 2a
\leq \kappa^{-1}I,\\
&\psi^{|\alpha|}(|D^{\alpha}a|+|D^{\alpha}b|+\psi |D^{\alpha}c|
+|D^{\alpha}\sigma|+\psi|D^{\alpha}\nu|)\leq K,\\
&\psi^{|\alpha|+1}D^{\alpha}F \in \mathbb{H}^{0}(\mathcal{D}),
\quad \psi^{|\beta|}D^{\beta}\phi\in L^{2}(\Omega, \mathscr{F}_{T},
L^{2}(\mathcal{D})),
\end{split}
\end{eqnarray*}
The first inequality is so called the super-parabolic condition (see e.g.
\cite[p139]{MaYo99}).
\end{ass}

Assumption \ref{ass:b3} implies that $a,b$ and $\sigma$ are bounded (by $K$).

\begin{ass}\label{ass:b4}
The functions $c,\nu$ and the derivatives of $a$ and $\sigma$ are bounded,
that is
\begin{equation*}
  |Da|+|D\sigma|+|c|+|\nu|\leq K.
\end{equation*}
\end{ass}

We will deal with several norms of various functions.
Denote $Q=Q(T)=(0,T) \times \mathcal{D}$.
For an integer $m$ and for real or complex-valued functions
$\phi,u$ defined on $\Omega \times \mathcal{D},
\Omega \times Q$ respectively we denote
\begin{eqnarray*}\begin{split}
&\|\phi\|_{\mathcal{D}}^{2} = \int_{\mathcal{D}} |\phi(x)|^{2}dx,&&
\|u\|_{Q}^{2} = \int_{Q} |u(t,x)|^{2}dx dt,\\
&[| \phi| ]_{m,\mathcal{D}}^{2} = \sum_{|\alpha|=m}
\int_{\mathcal{D}} |\psi^{m}D^{\alpha}\phi|^{2}dx,&&
[|u|]_{m,Q}^{2} = \sum_{|\alpha|=m}
\int_{Q} |\psi^{m}D^{\alpha}u|^{2}dx dt,\\
&\|\phi \|_{m,\mathcal{D}}^{2} = \sum_{r\leq m}[|\phi|]_{r,\mathcal{D}}^{2},&&
\|u \|_{m,Q}^{2} = \sum_{r\leq m}[|u|]_{r,Q}^{2},\\
&[\|\phi\|]_{m,\mathcal{D}}^{2} = E [|\phi|]_{m,\mathcal{D}}^{2},&&
[\|u\|]_{m,Q}^{2} = E [|u|]_{m,Q}^{2},\\
&\interleave \phi \interleave _{m,\mathcal{D}}^{2}
= E \|\phi\|_{m,\mathcal{D}}^{2},&&
\interleave u \interleave _{m,Q}^{2}= E \|u\|_{m,Q}^{2}.
\end{split}\end{eqnarray*}
The same notations will be used for vector-valued and matrix-valued
functions, and in the latter case we denote
$|u|^{2}=\sum\limits_{i,k}|u^{ik}|^{2}$.

Our main result is the following

\begin{thm}\label{thm:b1}
Let Assumptions \ref{ass:b1}-\ref{ass:b4} be satisfied. Then the equation
\begin{equation}\label{eq:b1}
dp=-\big[a^{ij}D_{ij}p+b^iD_{i}p-cp+\sigma^{ik}D_{i}q^k+\nu^kq^k+F\big]dt
+q^kdW^k_t
\end{equation}
with the boundary condition \eqref{con:b} has a unique solution $(p,q)$ such
that
$$ p\in \mathbb{H}^{1}_{0}(\mathcal{D}),\quad \psi
D^{2}p,q,\psi Dq\in \mathbb{H}^{0}(\mathcal{D}). $$ Moreover, for this
solution pair and any multi-index $\beta$ such that $|\beta|\leq n+1$, we
have
\begin{eqnarray}\label{prp:b1}
\psi^{|\beta|}D^{\beta}(D p),\;\psi^{|\beta|}D^{\beta}q
\in \mathbb{H}^{0}(\mathcal{D}),\quad
\psi^{|\beta|}D^{\beta}p \in C([0,T],L^{2}(\mathcal{D})) \ (a.s.),
\end{eqnarray}
and moreover
\begin{eqnarray}\label{leq:b1}
\interleave Dp\interleave^{2}_{n+1,Q}+E\sup_{t\leq
T}\|p(t,\cdot)\|^{2}_{n+1,\mathcal{D}}+\interleave q \interleave^{2}_{n+1,Q}
\leq C(\interleave \psi F\interleave^{2}_{n,Q}+\interleave
\phi\interleave^{2}_{n+1,\mathcal{D}}),
\end{eqnarray}
where the constant $C$ depends only on the norm of $\psi$ in
$C^{n+2}(\bar{\mathcal{D}})$, on the parameters $n,K,\rho_{0},\kappa$ and
$T$.
\end{thm}

The proof will be given in the final section. Note that the form of
Eq.\eqref{eq:b1} is slightly different from \eqref{eq:a1}, which is beside
the point because we could rewrite \eqref{eq:b1} into the form of
\eqref{eq:a1} due to Assumption \ref{ass:b4}. Now we give several remarks.

\begin{rmk}
Theorem \ref{thm:b1} implies the interior regularity of the equation
\eqref{eq:b1}. Indeed, for any domain $\mathcal{D}'$ such that
$\bar{\mathcal{D}'}\subset \mathcal{D}$ and for any multi-index $\beta$ such
that $|\beta|\leq n+1$ we have
\begin{eqnarray*}
\| D^{\beta}p_{x}\|^{2}_{\mathbb{H}^{0}(\mathcal{D}')}+E\sup_{t\leq
T}\|D^{\beta}p(t,\cdot)\|^{2}_{L^2(\mathcal{D}')}+\| D^{\beta}q\|^{2}
_{\mathbb{H}^{0}(\mathcal{D}')}\\
\leq C d^{-2(n+1)}(\interleave\psi F\interleave^{2}_{n,Q}
+\interleave \phi\interleave^{2}_{n+1,\mathcal{D}})
\end{eqnarray*}
with the same constant $C$ in \eqref{leq:b1} and with
$d=\textrm{dist}(\mathcal{D}',\partial\mathcal{D})$. As a corollary, if $n$
is large enough, then $(p,q)$ becomes a classical solution (see e.g.
\cite[p140]{MaYo99}) for Eq.\eqref{eq:b1} as a result of the Sobolev
Imbedding Theorem.
\end{rmk}

\begin{rmk}
Although we assumed that $d_1$ is finite, all constants denoted by $C$
will be independent of $d_1$, which allows us to extend our results
to the most general case of equations like \eqref{eq:b1}
with any Hilbert-space valued Wiener process. For instance, $d_1$ is
infinite while $q$ is $l^{2}$-valued process.
\end{rmk}

\begin{rmk}
With the help of the theory of interpolation, the main result in this paper
could be generalized to functions with fractional derivatives, thus allowing
$n$ to be a fraction rather than an integer.
\end{rmk}

Set $\|h\|^{2}=\int_{0}^{\infty}|h(x)|^{2}dx$. The next two lemmas proved in
\cite{Kryl94} will be used in several occations of our paper.

\begin{lem}\label{lem:b1}
\begin{enumerate}\item[\emph{(a)}] If $h(x)$ is a function of class
$H^{r+1}_{loc}(\mathbb{R}_{+}),r\geq 0$,and $h(x)=0$ for large x, then for
any $m>-\frac{1}{2}$
$$\int_{0}^{\infty}|x^{m}h^{(r)}(x)|^{2}dx\leq
(m+\frac{1}{2})^{-2}\int_{0}^{\infty}|x^{m+1}h^{(r+1)}(x)|^{2}dx.$$
\item[\emph{(b)}] If $h(x)$ is a function of class $H^{r}_{loc}(\mathbb{R}_{+}),r\geq
0$,and $h(x)=0$ for large x, and $p\geq 0,m>r-\frac{1}{2}$, then
$$\|x^{m+p}h^{(r)}\|\sim \|x^{m}(x^{p}h)^{(r)}\|,$$
in the sense that each norm times a constant, depending only on
m,p,r, is greater than the other one.
\item[\emph{(c)}] If $h(x)$ is a function of class $H^{1}_{loc}(\mathbb{R}_{+}),r\geq
0$,and $h(x)=0$ for large x, for any numbers $n\geq p\geq 0,n\geq 1$
we have$$\|x^{p}h\|\leq C(\|x^{n}h\|+\|h_{x}\|)$$ where
$C=C(n)$.\end{enumerate}
\end{lem}

\begin{lem}\label{lem:b2}
Let $v$ be a function on $\mathcal{D}$ and $\psi$ be the function defined
above.

\emph{(a)} If $v\psi^{-1},Dv\in L^{2}(\mathcal{D})$, then $v\in
H^{1}_{0}(\mathcal{D})$.

\emph{(b)} If $v,\psi Dv\in L^{2}(\mathcal{D})$, then $\psi v\in
H^{1}_{0}(\mathcal{D})$.

\emph{(c)} If $v\in H^{1}_{0}(\mathcal{D})$, then $v\psi^{-1}\in
L^{2}(\mathcal{D})$.

\emph{(d)} If $\psi v\in L^{2}(\mathcal{D})$, then $v\in
H^{-1}(\mathcal{D})$.
\end{lem}

\section{The 1-dimensional equation}

The investigation of one-dimensional equations in the next two sections is
most essential in this paper. We denote
$$\mathbb{R}_{+}=(0,\infty),\quad Q=(0,T)\times \mathbb{R}_{+},\quad
\interleave \cdot \interleave = \interleave \cdot \interleave_{Q},\quad
\|\cdot\|=\|\cdot\|_{\mathbb{R}_+}$$ in the next two sections.

\begin{thm}\label{thm:c1}
Fix constants $\kappa\in (0,1)$ and $\lambda_{0}\in \mathbb{R}_{+}$. Let $n$
and $r$ be integers s.t. $0\leq r\leq n$. Let $(p,q)$ be a solution of the
class $\mathbb{H}^{1}_{0}(\mathbb{R}_{+})\otimes
\mathbb{H}^{0}(\mathbb{R}_{+})$ for the equation
\begin{eqnarray}\label{eq:c1}
\left\{\begin{array}{l} dp=-[a^{2}p_{xx}+2\varepsilon a\lambda
p_{x}-\lambda^{2}p+\delta^{k}aq^{k}_{x}+\gamma^{k}\lambda
q^{k}+G_{x}]dt+q^{k}dW^{k}_{t}, \\
p(t,0)=0,\ p(T,x)=\phi(x),\quad t\in [0,T],\ x\in \mathbb{R}_{+},
\end{array}\right.
\end{eqnarray}
with the assumption that the predictable functions
$a,\lambda,\varepsilon,\delta$ and $\gamma$ take values in
$\mathbb{R},\mathbb{R},\mathbb{C},\mathbb{R}^{d_1}$ and $\mathbb{C}^{d_1}$
respectively and are bounded and all are independent of $x$. Suppose that
\begin{eqnarray}\label{con:c1}
\begin{split}
&\kappa \leq a \leq \kappa^{-1},\ \kappa \lambda_{0}\leq
\lambda \leq\kappa^{-1}\lambda_{0},\ |\delta|^{2}\leq 2-\kappa,\\
&x^{r}G^{(r)},\dots,G\in \mathbb{H}^{0}(\mathbb{R}_{+}),\quad
x^{r}\phi^{(r)},\dots,\phi\in
L^{2}(\Omega,\mathscr{F}_{T},L^{2}(\mathbb{R}_{+})),
\end{split}
\end{eqnarray}
and for any $z_{1},z_{2}\in \mathbb{C}$ and any $z_{3}\in \mathbb{C}^{d_1}$,
\begin{eqnarray}\label{con:c2}
\nonumber 2|z_{1}|^{2}+2|z_{2}|^{2}+|z_{3}|^{2}-4Re(\varepsilon
\bar{z}_{1}z_{2})
-2Re(\delta \bar{z}_{1}z_{3})-2Re(\gamma \bar{z}_{2}z_{3})\\
\geq \mu(\kappa)(|z_{1}|^{2}+|z_{2}|^{2}+|z_{3}|^{2})
\end{eqnarray}
with the constant $\mu(\kappa)>0$. Suppose finally that for a number $\rho
>0$ the functions $G(t,x),p(t,x)$ and $q(t,x)$ vanish if $x\geq \rho$. Then
\begin{eqnarray*}
\begin{split}
x^{r}p^{(r+1)},\dots,p_{x},x^{r}q^{(r)},\dots,q \in
\mathbb{H}^{0}(\mathbb{R}_{+}), \quad  x^{n}p^{(r)}\in
C([0,T],L^{2}(\mathbb{R}_{+}))\ (a.s.),
\end{split}
\end{eqnarray*}
and for a constant $C$ depending only on $n,r,\kappa$ but independent of
$\rho$ and $T$,
\begin{eqnarray}\label{leq:c1}
\begin{split}
\interleave x^{n}&p^{(r+1)}\interleave + \lambda_{0}\interleave
x^{n}p^{(r)}\interleave +\interleave x^{n}q^{(r)}\interleave
+(E\sup_{t\leq T}\|x^{n}p^{(r)}(t,\cdot)\|^{2})^{1/2}\\
\leq &C\big[(\lambda_{0}^{-1}\interleave x^{n}G^{(r+1)}\interleave)
\wedge\interleave x^{n}G^{(r)}\interleave+
(E\|x^{n}\phi^{(r)}\|^{2})^{1/2}\big]\\
&+C\lambda_{0}^{r-n}\big[(\lambda_{0}^{-1} \interleave
x^{r}G^{(r+1)}\interleave) \wedge\interleave x^{r}G^{(r)}\interleave+
(E\|x^{r}\phi^{(r)}\|^{2})^{1/2}\big],
\end{split}
\end{eqnarray}
where we put $\interleave x^{r}G^{(r+1)}\interleave =\infty$ if
$x^{r}G^{(r+1)} \notin \mathbb{H}^{0}(\mathbb{R}_{+})$.
\end{thm}

To prove this theorem, we need the following lemma, whose proof will be given
in the next section.

\begin{lem}\label{lem:c1}
Under the conditions of Theorem \ref{thm:c1}, we have
\begin{eqnarray*}
x^{r}p^{(r+1)},\dots,p_{x},\
x^{r}q^{(r)},\dots,q\in\mathbb{H}^{0}(\mathbb{R}_{+}).
\end{eqnarray*}
\end{lem}

\begin{proof}[Proof of Theorem \ref{thm:c1}]
Our proof consists of two steps. \medskip

\emph{Step 1}. Assume that $\lambda_{0}=1$.

First of all, from Lemma \ref{lem:c1} and from the fact
$$x^{j}p^{(j+2)}=(x^{j}p^{(j+1)})_{x}-j x^{j-1}p^{(j+1)},$$
it follows by induction that
$$x^{j}p^{(j+2)}\in H^{-1}(\mathbb{R}_{+}),\quad a.e.(\omega,t)\in \Omega\times[0,T],$$
for any $j=0,\dots,r$, and for any $\psi \in H^{1}_{0}(\mathbb{R}_{+})$
$$\int_{\mathbb{R}_{+}}\psi x^{j}p^{(j+2)} dx
= -\int_{\mathbb{R}_{+}}\psi_{x} x^{j}p^{(j+1)}dx -\int_{\mathbb{R}_{+}}\psi
j x^{j-1}p^{(j+1)} dx.$$ Obviously the function $x^{m}p^{(j+2)}\in
\mathbb{H}^{-1}(\mathbb{R}_{+})$ for $m\geq j$ (note that $p(t,x)=0$ if
$x>\rho$). Some similar arguments yield that
$$x^{m}G_{x}^{(j)}\in H^{-1}(\mathbb{R}_{+}),\quad a.e.(\omega,t)\in \Omega\times[0,T],$$
Moreover, in this situation for $a.e. (\omega,t)\in \Omega\times[0,T]$,
$$x^{m-1}p^{(j)}\in L^{2}(\mathbb{R}_{+}),\quad
(x^{m}p^{(j)})_{x} = m(x^{m-1}p^{(j)}) +(x^{m}p^{(j+1)})\in
L^{2}(\mathbb{R}_{+}),$$ which along with Lemma \ref{lem:b3}(a) implies
$$ x^{m}p^{(j)}\in H^{1}_{0}(\mathbb{R}_{+}),\quad a.e.(\omega,t)\in \Omega\times[0,T],$$

Consider the following equation
\begin{eqnarray*}
\begin{split}
d\big(x^{m}p^{(j)}\big)= & -\big{(}a^{2}x^{m}p_{xx}^{(j)}+2\varepsilon a
\lambda x^{m}p^{(j)}_{x}-\lambda^{2}x^{m}p^{(j)}\\
&+\delta a x^{m}q_{x}^{(j)}+\gamma \lambda
x^{m}q^{(j)}+x^{m}G_{x}^{(j)} \big{)}dt + x^{m}q^{(j)}dW_{t}.
\end{split}
\end{eqnarray*}
From the above arguments and Lemma \ref{lem:b5}, it follows that
$x^{m}p^{(j)}\in C([0,T ],L^{2}(\mathbb{R}_{+}))$ (a.s.). Applying It\^{o}'s
formula to evaluate $\|x^{m}p^{(j)}\|^{2}$, we get that for any $j\leq
m\wedge r$
\begin{eqnarray*}
\begin{split}
d\|x^{m}&p^{(j)}(t,\cdot)\|^{2}\\
=\ &\int_{\mathbb{R}_{+}}\bigg{[}2a^{2}|x^{m}p_{x}^{(j)}|^{2}
+4ma^{2}x^{2m-1}Re(\bar{p}^{(j)}p_{x}^{(j)})-4a\lambda x^{2m}Re(\varepsilon
\bar{p}^{(j)}p_{x}^{(j)})\\
&+2\lambda^{2}|x^{m}p^{(j)}|^{2} -2x^{2m}Re(\delta
a\bar{p}^{(j)}q_{x}^{(j)}+\gamma \lambda
\bar{p}^{(j)}q^{(j)})-2x^{2m}Re(\bar{p}^{(j)}G_{x}^{(j)})\\
&+|x^{m}q^{(j)}|^{2} \bigg{]}dx \cdot dt +
\int_{\mathbb{R}_{+}}2x^{2m}Re(\bar{p}^{(j)}q^{(j)})dx \cdot dW_{t}.
\end{split}
\end{eqnarray*}
Using the integration by parts, we have
\begin{eqnarray*}
\begin{split}
\int_{\mathbb{R}_{+}}4ma^{2}x^{2m-1}Re(\bar{p}^{(j)}p_{x}^{(j)})dx &=\
2a^{2}m(2m-1)\int_{\mathbb{R}_{+}}x^{2m-1}(|p^{(j)}|^{2})_{x} dx\\
&= - 2a^{2}m(2m-1)\|x^{m-1}p^{(j)}(t,\cdot)\|^{2}.
\end{split}
\end{eqnarray*}
Hence we get that
\begin{eqnarray}\label{eq:c2}
\begin{split}
-d\|x^{m}&p^{(j)}(t,\cdot)\|^{2}+\big{[}2a^{2}\|x^{m}p^{(j+1)}\|^{2}
+2\lambda^{2} \|x^{m}p^{(j)}\|^{2} + \|x^{m}q^{(j)}\|^{2} \big{]}dt\\
=\
&2a^{2}m(2m-1)\|x^{m-1}p^{(j)}\|^{2}dt+\int_{\mathbb{R}_{+}}\bigg{[}4a\lambda
x^{2m}Re(\varepsilon \bar{p}^{(j)}p_{x}^{(j)})\\
& + 2x^{2m}Re(\delta a\bar{p}_{x}^{(j)}q^{(j)} + \gamma \lambda
\bar{p}^{(j)}q^{(j)})
+ 4m x^{2m-1}Re(\delta a\bar{p}^{(j)}q^{(j)})\\
& +2x^{2m}Re(\bar{p}^{(j)}G_{x}^{(j)})\bigg{]}dx \cdot dt -
\int_{\mathbb{R}_{+}}2x^{2m}Re(\bar{p}^{(j)}q^{(j)})dx \cdot dW_{t}.
\end{split}
\end{eqnarray}
Applying the condition \eqref{con:c2} to \eqref{eq:c2} with
$z_{1}=ax^{m}p^{(j+1)},z_{2}=\lambda x^{m}p^{(j)},z_{3}=x^{m}q^{(j)}$,
we have
\begin{eqnarray*}
\begin{split}
-d\|x^{m}&p^{(j)}(t,\cdot)\|^{2}+\mu(\kappa)\big{[}a^{2}\|x^{m}p^{(j+1)}\|^{2}
+\lambda^{2} \|x^{m}p^{(j)}\|^{2} + \|x^{m}q^{(j)}\|^{2} \big{]}dt\\
\leq&\ 2a^{2}m(2m-1)\|x^{m-1}p^{(j)}\|^{2}dt+\int_{\mathbb{R}_{+}}\bigg{[}
4m x^{2m-1}Re(\delta a\bar{p}^{(j)}q^{(j)})\\
& +2x^{2m}Re(\bar{p}^{(j)}G_{x}^{(j)})\bigg{]}dx\cdot dt -
\int_{\mathbb{R}_{+}}2x^{2m}Re(\bar{p}^{(j)}q^{(j)})dx \cdot dW_{t}.
\end{split}
\end{eqnarray*}
By the Cauchy-Schwarz inequality, we get
\begin{eqnarray*}
\begin{split}
\int_{\mathbb{R}_{+}}4mx&^{2m-1}Re(\delta a\bar{p}^{(j)}q^{(j)})dx\\
\leq\ &\eta_{1}\cdot \delta^{2}\|x^{m}q^{(j)}(t,\cdot)\|^{2} +
\eta_{1}^{-1}\cdot 4m^{2}a^{2}\|x^{m-1}p^{(j)}(t,\cdot)\|^{2}\\
\int_{\mathbb{R}_{+}}2x^{2m}&Re(\bar{p}^{(j)}G_{x}^{(j)}) dx\\
\leq\ &\eta_{2}\cdot m^{2}\|x^{m-1}p^{(j)}(t,\cdot)\|^{2} +
\eta_{2}\|x^{m}p^{(j+1)}(t,\cdot)\|^{2} + 2\eta_{2}^{-1}\|x^{m}G^{(j)}(t,\cdot)\|^{2},\\
\int_{\mathbb{R}_{+}}2x^{2m}&Re(\bar{p}^{(j)}G_{x}^{(j)}) dx\\
\leq\ &\eta_{3}\|x^{m}p^{(j)}(t,\cdot)\|^{2} +
\eta_{3}^{-1}\|x^{m}G^{(j+1)}(t,\cdot)\|^{2}.
\end{split}
\end{eqnarray*}
Taking $\eta_{1},\eta_{2},\eta_{3}$ small enough and recalling
$\lambda_{0}=1$ in this step, we get
\begin{eqnarray}\label{leq:c2}
\begin{split}
-d\|x^{m}&p^{(j)}(t,\cdot)\|^{2}+\eta \big{[}\|x^{m}p^{(j+1)}\|^{2}
+\|x^{m}p^{(j)}\|^{2} + \|x^{m}q^{(j)}\|^{2} \big{]}dt\\
\leq\ &C\big{[}\|x^{m-1}p^{(j)}\|^{2}+\|x^{m}G^{(j)}\|^{2}\wedge
\|x^{m}G^{(j+1)}\|^{2} \big{]}dt + dM_{t},
\end{split}
\end{eqnarray}
where $M_{t}$ is a local martingale such that $dM_{t}$ is equal to the last
term in \eqref{eq:c2} and $\eta$ depends only on $\kappa$. In fact it is not
hard to prove that $M_{t}$ is a uniformly integrable martingale by the
Burkholder-Davis-Gundy inequality. Then integrating this inequality with
respect to $t$ and taking expectations, we obtain
\begin{equation}\label{leq:c3}
T_{m}^{j} \leq C T_{m-1}^{j-1} + C R_{m}^{j},
\end{equation}
where
\begin{eqnarray*}
\begin{split}
&T_{m}^{j}= \interleave x^{m}p^{(j+1)} \interleave^{2}
+ \interleave x^{m}p^{(j)}\interleave^{2}
+ \interleave x^{m}q^{(j)}\interleave^{2},\\
&R_{m}^{j}= \interleave x^{m}G^{(j)}\interleave^{2}\wedge\interleave
x^{m}G^{(j+1)} \interleave^{2} + E\| x^{m}\phi^{(j)}\|^{2}.
\end{split}
\end{eqnarray*}
In particular, it follows from the similar argument that
$$T_{0}^{0}\leq C R^{0}_{0},$$
which along with \eqref{leq:c3} yields
\begin{eqnarray*}
\begin{split}
T_{m}^{j} &\leq C T_{m-j}^{0} + C (
R_{m}^{j}+\cdots + R_{0}^{m-j}) \\
&\leq C (R_{m-j}^{0}+\cdots +R_{0}^{0}) + C (R_{m}^{j}+\cdots + R_{0}^{m-j}),
\end{split}
\end{eqnarray*}
where the constant $C$ depends
only on $m,r,\kappa$ but is independent of $\rho$ and $T$. In particular for
$j=r,m=n$ and by Lemma \ref{lem:b1}(a) we have
\begin{eqnarray}\label{leq:c4}
\begin{split}
T_{n}^{r}\ &\leq\ C (R_{n-r}^{0}+\cdots +R_{0}^{0}) + C (
R_{n}^{r}+\cdots + R_{0}^{n-r})\\
&\leq\ C(R_{n-r}^{0}+\cdots +R_{0}^{0}) + C R_{n}^{r}.
\end{split}
\end{eqnarray}

In order to estimate the last term on the left in \eqref{leq:c1} we
come back to \eqref{leq:c2} with $j=r,m=n$, by the
Burkholder-Davis-Gundy inequality we have
\begin{equation}\label{leq:c5}
E\sup_{t\leq T}\|x^{n}p^{(r)}(t,\cdot)\|^{2}\ \leq\
CT_{n-1}^{r-1}+CR_{n}^{r}+CS_{n}^{r},
\end{equation}
where
\begin{eqnarray*}
\begin{split}
S_{n}^{r}\ =&\
E\bigg{[}\int_{0}^{T}\bigg{(}\int_{R_{+}}2x^{2n}Re(\bar{p}^{(r)}q^{(r)})dx
\bigg{)}^{2}dt \bigg{]}^{1/2}\\
\leq&\ C\big{(}\interleave x^{n}p^{(r)} \interleave^{2} +
\interleave
x^{n}q^{(r)}\interleave^{2} \big{)}\\
\leq&\ CT_{n}^{r}.
\end{split}
\end{eqnarray*}
Combining \eqref{leq:c4} and \eqref{leq:c5}, we have
$$T_{n}^{r}+E\sup_{t\leq T}\|x^{n}p^{(r)}(t,\cdot)\|^{2}
\leq C(R_{n-r}^{0}+\cdots +R_{0}^{0}) + C R_{n}^{r},$$ where $C=C(n,\kappa)$.

Now we have to get rid of the term $R_{n-r}^{0}+\cdots +R_{0}^{0}$. From
Lemma \ref{lem:b1}, it follows that $R_{n-r}^{0}\leq CR_{n}^{r}$, and
moreover for $0\leq m \leq n-r-1$,
\begin{eqnarray*}
\begin{split}
\|x^{m}\phi\| &\leq C\|x^{m+r}\phi^{(r)}\| \leq
C(\|x^{n}\phi^{(r)}\|+\|x^{r}\phi^{(r)}\|),\\
\interleave x^{m}G\interleave &\leq C (\interleave
x^{n}G^{(r)}\interleave + \interleave x^{r}G^{(r)} \interleave),\\
\interleave x^{m}G\interleave &\leq C \interleave
x^{m+r}G^{(r)}\interleave \leq C (\interleave
x^{n-1}G^{(r)}\interleave + \interleave x^{r}G^{(r)} \interleave)\\
&\leq C (\interleave
x^{n}G^{(r+1)}\interleave + \interleave x^{r}G^{(r)} \interleave),\\
\interleave x^{m}G\interleave &\leq C \interleave
x^{m+r+1}G^{(r+1)}\interleave \leq C (\interleave
x^{n}G^{(r+1)}\interleave + \interleave x^{r}G^{(r+1)} \interleave),\\
\interleave x^{m}G\interleave &\leq C (\interleave x^{n-r}G
\interleave + \interleave G_{x} \interleave) \leq C (\interleave
x^{n}G^{(r)}\interleave + \interleave x^{r}G^{(r+1)} \interleave).\\
\end{split}
\end{eqnarray*}
Therefore
\begin{eqnarray*}
\begin{split}
&\interleave x^{m}G\interleave \wedge \interleave
x^{m}G_{x}\interleave\ \leq\ \interleave x^{m}G\interleave\\
\leq\ & C \big{[}(\interleave x^{n}G^{(r+1)}\interleave \wedge
\interleave x^{n}G^{(r)}\interleave)+(\interleave
x^{r}G^{(r+1)}\interleave \wedge \interleave
x^{r}G^{(r)}\interleave) \big{]}.
\end{split}
\end{eqnarray*}
So the inequality \eqref{leq:c1} is proved in this situation. \medskip

\emph{Step 2}. The general case can be reduced to the case in Step 1 by
scaling, that is by introducing the new functions
$$
\left.\begin{array}{ll}
\tilde{p}(t,x)=p(\lambda_{0}^{-2}t,\lambda_{0}^{-1}x),\quad &\tilde{q}(t,x)=
\lambda_{0}^{-1} q(\lambda_{0}^{-2}t,\lambda_{0}^{-1}x),\\
\tilde{G}(t,x)=\lambda_{0}^{-1}G(\lambda_{0}^{-2}t,\lambda_{0}^{-1}x),\quad
&\tilde{\phi}(x)=\phi (\lambda_{0}^{-1}x),
\end{array}\right.
$$
and a new Wiener process
$\tilde{W}_{t}=\lambda_{0}W_{\lambda_{0}^{-2}t}$. It is not hard to
show that $\tilde{p},\tilde{q}$ satisfies the equation
\begin{equation*}
\left\{\begin{array}{ll} d\tilde{p}\ =\ -\big{[}a^{2}\tilde{p}_{xx} +
2\varepsilon a(\lambda\lambda_{0}^{-1})\tilde{p}_{x} -
(\lambda\lambda_{0}^{-1})^{2}\tilde{p}\\
\qquad\qquad + \delta a\tilde{q}_{x} + \gamma
(\lambda\lambda_{0}^{-1})\tilde{q} + \tilde{G}_{x} \big{]}dt +
\tilde{q} d\tilde{W}_{t},\\
\tilde{p}(t,0)=0,\quad
\tilde{p}(\lambda_{0}^{2}T,x)=\tilde{\phi}(x).
\end{array}\right.
\end{equation*}
Direct calculus shows that
\begin{eqnarray*}
\left.\begin{array}{ll}
\interleave x^{n}\tilde{p}^{(r)}\interleave =
\lambda_{0}^{n-r+3/2}\interleave x^{n}p^{(r)} \interleave,
&\interleave x^{n}\tilde{q}^{(r)}\interleave =
\lambda_{0}^{n-r+1/2}\interleave x^{n}q^{(r)} \interleave,\\
\interleave x^{n}\tilde{G}^{(r)}\interleave =
\lambda_{0}^{n-r+1/2}\interleave x^{n}G^{(r)} \interleave,
&\| x^{n}\tilde{\phi}^{(r)}\| = \lambda_{0}^{n-r+1/2}\| x^{n}\phi^{(r)} \|.
\end{array}\right.
\end{eqnarray*}
Note that the constant $C$ in Step 1 is independent of the terminal time $T$
and $\rho$. Hence the inequality \eqref{leq:c1} is easily obtained from the
result in Step 1. The proof is complete.
\end{proof}

\begin{rmk}
With the help of Lemma \ref{lem:b1}(a),(b), we can rewrite \eqref{lem:c1}
into some convenient forms which will be used in Section 5 directly.
Replacing $n$ and $r$ in \eqref{leq:c1} by $n+1$ and $r+1$ respectively, we
get
\begin{eqnarray}\label{leq:c6}
\begin{split}
&\interleave x^{n}(x p_{xx})^{(r)}\interleave + \lambda_{0}\interleave
x^{n}(x p_{x})^{(r)}\interleave \\&+\interleave x^{n}(x
q_{x})^{(r)}\interleave
+(E\sup_{t\leq T}\|x^{n}(x p_{x})^{(r)}(t,\cdot)\|^{2})^{1/2}\\
&~~\leq C\big[\interleave x^{n}(x G_{x})^{(r)}\interleave+
(E\|x^{n}(x\phi_{x})^{(r)}\|^{2})^{1/2}\big]\\
&~~~~~~+C\lambda_{0}^{r-n}\big[ \interleave x^{n}(x G_{x})^{(n)}\interleave +
(E\|x^{n}(x\phi_{x})^{(n)}\|^{2})^{1/2}\big],
\end{split}
\end{eqnarray}
for $0\leq r\leq n$, if the right-hand side is finite. Furthermore, replacing
$n$ in \eqref{leq:c1} by $n+1$, we get
\begin{eqnarray}\label{leq:c7}
\begin{split}
&\interleave x^{n}(x p_{x})^{(r)}\interleave + \lambda_{0}\interleave x^{n}(x
p)^{(r)}\interleave \\&+\interleave x^{n}(x q)^{(r)}\interleave
+(E\sup_{t\leq T}\|x^{n}(x p)^{(r)}(t,\cdot)\|^{2})^{1/2}\\
&~~\leq C\big[\lambda_{0}^{-1}\interleave x^{n}(x G_{x})^{(r)}\interleave +
(E\|x^{n}(x\phi)^{(r)}\|^{2})^{1/2}\big]\\
&~~~~~~+C\lambda_{0}^{r-n-1}\big[\interleave x^{n}(x G_{x})^{(n)}\interleave+
(E\|x^{n}(x\phi_{x})^{(n)}\|^{2})^{1/2}\big],
\end{split}
\end{eqnarray}
for $0\leq r\leq n$, if the right-hand side is finite.
\end{rmk}

\section{Proof of Lemma \ref{lem:c1}}

We need several lemmas.
First we define the difference quotient (see e.g. \cite{Evan98}) by
$$\nabla _{h}u(x) = \frac{u(x+h)-u(x-h)}{2h},\quad h\neq 0$$
for $u\in L^{2}(\mathbb{R})$. The following basic lemma about the difference
quotient can be found in any fundamental textbooks on PDEs (e.g.
\cite{Evan98}).
\begin{lem}\label{lem:d1}
If $u\in H^{1}(\mathbb{R})$, then $\|\nabla _{h}u\|\leq \|u_{x}\|$, and
$\nabla _{h}u\rightarrow u_{x}$ strongly in $L^{2}(\mathbb{R})$ as
$h\downarrow 0$.
\end{lem}

To deal with the component $q$, we need the following

\begin{lem}\label{lem:d2}
If $u\in L^{2}(\mathbb{R})$, then $\|\nabla _{h}u\|_{-1}\leq \|u\|$, and
$\nabla _{h}u\rightarrow u_{x}$ strongly in $H^{-1}(\mathbb{R})$ as
$h\downarrow 0$.
\end{lem}

\begin{proof}
The first assertion follows from the previous lemma by duality. To
prove the second assertion, we only need to verify that the Bessel
potential $(1-\Delta)^{-1/2}$ is commutative with the difference
quotient $\nabla _{h}$ and the differential operator $D_{x}$. It obviously
holds true since that the Fourier transforms of the three operator are
$$(1+|\xi|^{2})^{-1/2},~~ \frac{e^{i\xi h}-e^{-i\xi h}}{2h},~~ i\xi$$
respectively which are multipliers.
\end{proof}

Applying the Lebesgue's dominated convergence theorem, we easily obtain

\begin{cor}\label{cor:d1}
\begin{enumerate}\item[\emph{(a)}]
If $p\in \mathbb{H}^{1}_{0}(\mathbb{R})$, then $$\| \nabla _{h}p\|
_{\mathbb{H}^{0}(\mathbb{R})}\leq \| p_{x}\| _{\mathbb{H}^{0}(\mathbb{R})},$$
and $\nabla _{h}p\rightarrow p_{x}$ strongly in $\mathbb{H}^{0}(\mathbb{R})$
as $h\downarrow 0$.
\item[\emph{(b)}] If $q\in \mathbb{H}^{0}(\mathbb{R})$, then
$$\| \nabla _{h}q\| _{\mathbb{H}^{-1}(\mathbb{R})}\leq \| q\|
_{\mathbb{H}^{0}(\mathbb{R})},$$ and $\nabla _{h}q\rightarrow q_{x}$ strongly
in $\mathbb{H}^{-1}(\mathbb{R})$ as $h\downarrow 0$. \end{enumerate}
\end{cor}

\begin{rmk}\label{rem:d1}
If the functions $p$ and $q$ in Corollary \ref{cor:d1} vanish for $|x|\geq
\rho$, where $\rho$ is a large number, it is easy to check that $x \nabla
_{h}p\rightarrow xp_{x}$ strongly in $\mathbb{H}^{0}(\mathbb{R})$ and
$x\nabla _{h}q\rightarrow xq_{x}$ strongly in $\mathbb{H}^{-1}(\mathbb{R})$
as $h\downarrow 0$.
\end{rmk}

The next lemma is known in the theory of Sobolev spaces. Note that a function
in $H^{1}(\mathbb{R})$ belongs to $C^{\frac{1}{2}}(\mathbb{R})$.

\begin{lem}\label{lem:b4}
Let $u\in H^{1}(\mathbb{R})$ be an odd function, then its restriction on
$\mathbb{R}_{+}$ belongs to $H^{1}_{0}(\mathbb{R}_+)$.
\end{lem}

We also need the following lemma concerning the odd (or even) continuation of
a function, in which we denote $\bar{f}$ as the odd (or even) continuation of
the function $f$ defined on $\mathbb{R}_+$.

\begin{lem}\label{lem:d3}
Let $n$ be an integer. Then $$f,\dots,x^{n}f^{(n)}\in L^{2}(\mathbb{R}_{+})
~~ \Leftrightarrow ~~ \bar{f},\dots,x^{n}\bar{f}^{(n)}\in
L^{2}(\mathbb{R}).$$
\end{lem}

\begin{proof}
It is a direct result from Lemma \ref{lem:b2}(b). Indeed, without loss of
generality, we assume that $f$ vanishes for large $x$. Then $f, x f'\in
L^{2}(\mathbb{R}_{+})$ implies $ x f \in H^{1}_{0}(\mathbb{R}_+)$, and then
$x \bar{f} \in H^{1}(\mathbb{R})$ owning to the theory of Sobolev spaces.
Thus $x\bar{f}'\in L^{2}(\mathbb{R})$. Then the necessity easily follows from
induction. The sufficiency is obvious.
\end{proof}

Furthermore, it is easy to show that
$x^{r}\bar{f}^{(r)}=\overline{x^{r}f^{(r)}}$ for any $r\leq n$, where
$\overline{x^{r}f^{(r)}}$ is the odd (or even) continuation of
$x^{r}f^{(r)}$. Hence, for any $r\leq n$ we have
$$\|x^{r}\bar{f}^{(r)}\|_{L^2(\mathbb{R})} = 2\|x^{r}f^{(r)}\|_{L^2(\mathbb{R}_+)}.$$

\begin{proof}[Proof of Lemma \ref{lem:c1}]
The proof consists of the following three steps. \medskip

\emph{Step 1}. Assume that $a=\lambda,\varepsilon=\beta=0$ and $\gamma=0$. In
this step, we denote
$$\interleave\cdot\interleave_{\mathbb{R}}=
\interleave\cdot\interleave_{(0,T)\times\mathbb{R}}.$$

First from Lemma \ref{lem:b3} we see that $(p,q)$ is the unique (weak)
solution of Eq.\eqref{eq:c1} in the space
$\mathbb{H}^{1}_{0}(\mathbb{\mathbb{R}}_+)\otimes\mathbb{H}^{0}
(\mathbb{\mathbb{R}}_+)$.

Let $\bar{G}$ be the even continuation of $G$ and
$\bar{\phi}$ be the odd continuation of $\phi$. It follows from
Lemma \ref{lem:d3} that
\begin{eqnarray}
\bar{G},\dots,x^{n}\bar{G}^{(n)} \in \mathbb{H}^{0}(\mathbb{R}),\quad
\bar{\phi} ,\dots,x^{r}\bar{\phi}^{(r)} \in
L^{2}(\Omega,\mathscr{F}_{T},L^{2}(\mathbb{\mathbb{R}}_{+})).
\end{eqnarray}
Thus according to Lemma \ref{lem:b3} the equation
\begin{eqnarray}\label{eq:d1}
\left\{\begin{array}{l}
d\bar{p}=-[a^{2}\bar{p}_{xx}-a^{2}\bar{p}+\delta^{k}a\bar{q}^{k}_{x}
+\bar{G}_{x}]dt+\bar{q}^{k}dW^{k}_{t}, \\
\bar{p}(T,x)=\bar{\phi}(x),\quad t\in [0,T],\ x\in \mathbb{R},
\end{array}\right.
\end{eqnarray}
has a unique solution pair $$(\bar{p},\bar{q})\in
\mathbb{H}^{1}_{0}(\mathbb{R})\otimes \mathbb{H}^{0}(\mathbb{R}),$$ such that
\begin{eqnarray}\label{leq:d1}
\begin{split} \interleave \bar{p}_{x} \interleave_{\mathbb{R}}^{2}+ E\sup_{t\leq
T}\|\bar{p}(t,\cdot)\|_{\mathbb{R}}^{2}+\interleave
\bar{q}\interleave_{\mathbb{R}}^{2} \leq &C\big{(}\|
\bar{G}_{x}\|_{\mathbb{H}^{-1}(\mathbb{R})}^{2} + E\|
\bar{\phi}\|_{\mathbb{R}}^{2}\big{)},\\
\leq &C\big{(}\interleave \bar{G}\interleave_{\mathbb{R}}^{2} + E\|
\bar{\phi}\|_{\mathbb{R}}^{2}\big{)}. \end{split}
\end{eqnarray}
with the constant $C=C(\kappa,T)$.

By symmetry the above functions $\bar{p}$ and $\bar{q}$ are odd with respect
to $x$, and then from Lemma \ref{lem:b4}, their restrictions on
$\mathbb{\mathbb{R}}_{+}$ belong to
$\mathbb{H}^{1}_{0}(\mathbb{\mathbb{R}}_+)\otimes
\mathbb{H}^{0}(\mathbb{\mathbb{R}}_+)$, which implies that
$\bar{p}=p,\bar{q}=q$ on $\mathbb{\mathbb{R}}_+$ (by the uniqueness of the
weak solution).

Now applying $x\nabla_{h}$ on both sides of Eq.\eqref{eq:d1}, we have
\begin{equation*}
\left\{\begin{array}{l}
d(x\nabla_{h}\bar{p})=-\big{[}a^{2}(x\nabla_{h}\bar{p})_{xx} -
a^{2} x\nabla_{h}\bar{p}+\delta a(x\nabla_{h}\bar{q})_{x}\\
\qquad\qquad\quad -2a^{2}(\nabla_{h}\bar{p})_{x}-\delta a \nabla_{h}\bar{q}
+(x\nabla_{h}\bar{G})_{x}-\nabla_{h}\bar{G}\big{]}dt
+x\nabla_{h}\bar{q}dW_{t},\\
(x\nabla_{h}\bar{p})|_{t=T}=x\nabla_{h}\bar{\phi}.
\end{array}\right.
\end{equation*}
In view of Lemma \ref{lem:d2}, it follows from Corollary \ref{cor:d1} and
\eqref{leq:d1} that
\begin{eqnarray}\label{leq:d2}
\begin{split}
&\interleave (x\nabla_{h}\bar{p})_{x}
\interleave_{\mathbb{R}}^{2}+\interleave
x\nabla_{h}\bar{q}\interleave_{\mathbb{R}}^{2}\\
\leq\ & C\big{(}\interleave \nabla_{h}\bar{p}\interleave_{\mathbb{R}}^{2}+ \|
\nabla_{h}\bar{q} \|^{2}_{\mathbb{H}^{-1}(\mathbb{R})}+ \interleave
x\nabla_{h}\bar{G}\interleave_{\mathbb{R}}^{2} +\|
\nabla_{h}\bar{G}\|^{2}_{\mathbb{H}^{-1}(\mathbb{R})}+
E\| x\nabla_{h}\bar{\phi}\|_{\mathbb{R}}^{2}\big{)}\\
\leq\ & C\big{(}\interleave \bar{p}_{x}\interleave_{\mathbb{R}}^{2}+
\interleave \bar{q}\interleave_{\mathbb{R}}^{2}+ \interleave
x\bar{G}_{x}\interleave_{\mathbb{R}}^{2} +\interleave
\bar{G}\interleave_{\mathbb{R}}^{2}+
E\| x\bar{\phi}_{x}\|_{\mathbb{R}}^{2}\big{)}\\
\leq\ & C\big{(}\interleave x\bar{G}_{x}\interleave_{\mathbb{R}}^{2}
+\interleave \bar{G}\interleave_{\mathbb{R}}^{2}+ E\|
x\bar{\phi}_{x}\|_{\mathbb{R}}^{2}+E\| \bar{\phi}\|_{\mathbb{R}}^{2}\big{)}.
\end{split}
\end{eqnarray}
Therefore, $\{x\nabla_{h}\bar{p}\}$ are bounded in
$\mathbb{H}^{1}_{0}(\mathbb{R})$ and $\{x\nabla_{h}\bar{q}\}$ are bounded in
$\mathbb{H}^{0}(\mathbb{R})$, uniformly with respect to $h$. Thus there exist
two functions
$$u\in \mathbb{H}^{1}_{0}(\mathbb{R}),\quad v\in \mathbb{H}^{0}(\mathbb{R}),$$
which are respectively the weak limits of (subsequences of)
$\{x\nabla_{h}\bar{p}\}$ and $\{x\nabla_{h}\bar{q}\}$ in as $h\downarrow 0$.

On the other hand, it follows from Remark \ref{rem:d1} that $x \nabla
_{h}\bar{p}\rightarrow x\bar{p}_{x}$ strongly in $\mathbb{H}^{0}(\mathbb{R})$
and $x\nabla _{h}\bar{q}\rightarrow x\bar{q}_{x}$ strongly in
$\mathbb{H}^{-1}(\mathbb{R})$ as $h\downarrow 0$. According to the uniqueness
of the limit, we have
$$ x\bar{p}_{x}=u  \in \mathbb{H}^{1}_{0}(\mathbb{R}),\quad
x\bar{q}_{x}=v  \in \mathbb{H}^{0}(\mathbb{R}), $$ which implies that
$$ xp_{xx}, xq_{x}\in \mathbb{H}^{0}(\mathbb{\mathbb{R}}_{+}).$$
Moveover, it follows from \eqref{leq:d2} that
\begin{equation}\label{leq:d3}
\interleave xp_{xx} \interleave^{2}+\interleave xq_{x}\interleave^{2} \leq\
C\big{(}\interleave xG_{x}\interleave^{2} +\interleave G\interleave^{2}+ E\|
x\phi_{x}\|^{2}+E\| \phi \|^{2}\big{)}.
\end{equation}

\emph{Step 2}. Now we prove by induction the assertion of Lemma \ref{lem:c1}
under the assumption that $a=\lambda,\varepsilon=\beta=0,\gamma=0$.

First the assertion is proved for $r=1$ in Step 1.

Assume that the assertion holds true for some $r\geq 1$, that is the
condition \eqref{con:c1} implies the following
\begin{equation}\label{prp:d1}
x^{r}p^{(r+1)},\dots,p_{x},x^{r}q^{(r)},\dots,q\in
\mathbb{H}^{0}(\mathbb{\mathbb{R}}_{+}).
\end{equation}
Since the odd function $x\bar{p}_{x}\in \mathbb{H}^{1}_{0}(\mathbb{R})$ (see
Step 1), we have $x p_{x}\in \mathbb{H}^{1}_{0}(\mathbb{\mathbb{R}}_+)$. It
is not hard to show the function pair $(x p_{x}, x q_{x})$ of the class
$\mathbb{H}^{1}_{0}(\mathbb{\mathbb{R}}_+)\otimes\mathbb{H}^{0}(\mathbb{\mathbb{R}}_+)$
solves the equation
\begin{equation}\label{eq:d2}
\left\{\begin{array}{l} d u=-\big{[}a^{2}u_{xx} - a^{2} u+\delta a v_{x}
-2a^{2}p_{xx}-\delta a q_{x} +x G_{xx}\big{]}dt
+ v dW_{t},\\
u|_{t=T}=x\phi_{x}
\end{array}\right.
\end{equation}
with unknown functions $u,v$. From the assumption \eqref{prp:d1} and the
assumption \eqref{con:c1} on $G,\phi$ and $r+1$, we see that the function
$$-2a^{2}p_{xx}-\delta a q_{x} +x G_{xx}=\big( -2a^{2}p_{x}-\delta a q_{x}
+ x G_{x}\big)_{x} - G_{x}$$ and the integer $r$ satisfy the condition
\eqref{con:c1}. Applying our assumption to Eq.\eqref{eq:d2}, we have
$$x^{r}(x p_{x})^{(r+1)},\dots,(x p_{x})_{x},x^{r}(xq_{x})^{(r)},\dots,xq_{x}
\in \mathbb{H}^{0}(\mathbb{\mathbb{R}}_{+}),$$ which implies the assertion of
Lemma \ref{lem:c1} as a result of Lemma \ref{lem:b1}(b). \medskip

\emph{Step 3}. Now consider the general situation. Rewrite
Eq.\eqref{eq:d1} as
\begin{eqnarray*}
\begin{split}
dp=-&\big{\{} a^{2}p_{xx}-a^{2}p+\delta^{k}aq^{k}_{x}+[2\varepsilon a \lambda p
+(\lambda^{2}-a^{2}+i\beta)P\\
&-\gamma^{k}\lambda Q +G]_{x} \big{\}} dt + q^{k}dW_{t}^{k},
\end{split}
\end{eqnarray*}
where
$$P(t,x)=\int_{x}^{\infty}p(t,y)dy,\quad Q(t,x)=\int_{x}^{\infty}q(t,y)dy.$$
Note that $P(t,x)=Q^{k}(t,x)=0$ for $x\geq \rho$ and
\begin{eqnarray*}
2\varepsilon a \lambda p +(\lambda^{2}-a^{2}+i\beta)P -\gamma^{k}\lambda Q +G
\in \mathbb{H}^{0}(\mathbb{R}_{+}),\\
x[2\varepsilon a \lambda p +(\lambda^{2}-a^{2}+i\beta)P -\gamma^{k}\lambda Q
+G]_{x} \in \mathbb{H}^{0}(\mathbb{R}_{+}).
\end{eqnarray*}
It follows from Step 1 that
$$ xp_{xx},p_{x},xq_{x},q \in \mathbb{H}^{0}(\mathbb{R}_{+}).$$

Identically as in Step 2, we prove by induction that Lemma \ref{lem:c1} holds
true in the general case for any $r$. The proof of the lemma is complete.
\end{proof}

\section{The equation with coefficients independent of $x$}

In this section we are concerned with the equation whose coefficients are
independent of $x$ in the half space. Denote $\mathbb{R}^{d}_{+}=\{x\in
\mathbb{R}^{d}:x^{1}>0\}$.

\begin{thm}\label{thm:e1}
Consider the following equation
\begin{equation}\label{eq:e1}
\left\{\begin{array}{l}
dp=-[a^{ij}D_{ij}p+\sigma^{ik}D_{i}q^k+F]dt+q^kdW^k_t, \\
p(t,x)=0,\ x\in\partial \mathbb{R}^{d}_{+},\\
p(T,x)=\phi(x),\ x\in \mathbb{R}^{d}_{+}.
\end{array}\right.
\end{equation}
Assume that $a$ and $\sigma$ satisfy Assumption \ref{ass:b1} and, in
addition, they are bounded and independent of $x$ and
\begin{equation}\label{con:e1} \kappa I+\sigma\sigma^{*}\leq 2a
\leq \kappa^{-1}I,\ \forall (\omega ,t).\end{equation} Suppose that for any
multi-indices $\alpha$ and $\beta$ such that $|\alpha|\leq n$ and
$|\beta|\leq n+1$ we have
\begin{eqnarray}\label{con:e2}
\tilde{\psi}^{|\alpha|+1}D^{\alpha}F \in
\mathbb{H}^{0}(\mathbb{R}^{d}_{+}),\quad \tilde{\psi}^{|\beta|}D^{\beta}\phi
\in L^{2}(\Omega, \mathscr{F}_{T},L^{2}(\mathbb{R}^{d}_{+})),
\end{eqnarray}
where $\tilde{\psi} (x)=x^{1}$. Then equation \eqref{eq:e1} has a unique
solution $(p,q)$ such that
$$p \in \mathbb{H}^{1}_{0}(\mathbb{R}^{d}_{+}),\quad \tilde{\psi} D^{2}p,q,\tilde{\psi} Dq
\in \mathbb{H}^{0}(\mathbb{R}^{d}_{+}).$$ For this solution and any
multi-index $\beta$ such that $|\beta|\leq n+1$, we have
\begin{eqnarray}\label{prp:e1}
\tilde{\psi}^{|\beta|}D^{\beta}(D p),\tilde{\psi}^{|\beta|}D^{\beta}q\in
\mathbb{H}^{0}(\mathbb{R}^{d}_{+}),\quad \tilde{\psi}^{|\beta|}D^{\beta}p\in
C([0,T],L^{2}(\mathbb{R}^{d}_{+})) \quad (a.s.),
\end{eqnarray}
and moreover
\begin{eqnarray}\label{leq:e}
[\|\tilde{\psi} D^{2}p\|]^{2}_{m,Q}+E\sup_{t\leq T}[\|\tilde{\psi} D
p(t,\cdot)\|]^{2}_{m,\mathbb{R}^{d}_{+}}+[\|
\tilde{\psi} D q\|]^{2}_{m,Q}\\
\nonumber \leq C([\| \tilde{\psi} F\|]^{2}_{m,Q}+\interleave \phi
\interleave^{2}_{m+1,\mathbb{R}^{d}_{+}}),
\end{eqnarray}
where $\tilde{\psi}(x)=x^{1},Q=(0,T)\times \mathbb{R}^{d}_{+},m=0,\dots,n$,
and the constant $C$ depends only on $d,n,K$ and $\kappa$.
\end{thm}

\begin{proof}
The proof consists of three steps.\medskip

\emph{Step 1}. We make the following additional assumption in this step.
\medskip

(A) \textit{Eq.\eqref{eq:e1} has a unique solution $(p,q)\in
\mathbb{H}^{1}_{0}(\mathbb{R}^{d}_{+})\otimes
\mathbb{H}^{0}(\mathbb{R}^{d}_{+})$, and for a number $\rho >0$ the functions
$F(t,x),p(t,x)$ and $q(t,x)$ vanish if $x^{1}\geq \rho$.} \medskip

Put
$$G(t,x)= -\int_{x^{1}}^{\infty} F(t,s,x^{2},\dots,x^{d}) ds,\quad
(t,x)\in [0,T]\times \mathbb{R}^{d}_{+},$$ which belongs to
$\mathbb{H}^{0}(\mathbb{R}^{d}_{+})$. Indeed, note that for almost every
$(\omega,t,x^{2},\dots,x^{d})$, $F$ as a function with respect to $x^{1}$
belongs to $L^{2}_{loc}(\mathbb{R}_{+})$, which implies
$G(\omega,t,\cdot,x^{2},\dots,x^{d}) \in H^{1}_{loc}(\mathbb{R}_{+})$. With
the help of Lemma \ref{lem:b1}(b) it is not hard to show that $\tilde{\psi}
G_{x^1}$ and $(\tilde{\psi} G)_{x^1}$ are of the same class while the former
belongs to $\mathbb{H}^{0}(\mathbb{R}^{d}_{+})$, which implies $G\in
\mathbb{H}^{0}(\mathbb{R}^{d}_{+})$. Moreover
$\tilde{\psi}^{|\beta|}D^{\beta}G\in \mathbb{H}^{0}(\mathbb{R}^{d}_{+})$ for
any $\beta$ such that $|\beta|\leq n+1$.

Now for a function $u(x)$ defined on $\mathbb{R}^{d}_{+}$, denote by
$\hat{u}(x^{1},\xi^{2},\dots,\xi^{d})$ its Fourier transform with respect to
$(x^{2},\dots,x^{d})$. Then we see that for almost every
$\xi=(\xi^{2},\dots,\xi^{d})\in \mathbb{R}^{d-1}\backslash\{0\}$, the
functions $\hat{p}=\hat{p}(t,x^{1},\xi^{2},\dots,\xi^{d})$ and
$\hat{q}=\hat{q}(t,x^{1},\xi^{2},\dots,\xi^{d})$ satisfy the equation
\begin{equation}\label{eq:e2}
\left\{\begin{array}{l}
d\hat{p}=-\big{[}a^{11}\hat{p}_{x^1x^1}+2\tilde{b}^{1}\hat{p}_{x^1}
-\tilde{c}\hat{p} +
\sigma^{1k}\hat{q}_{x^1}^{k}+\tilde{\nu}^{k}\hat{q}^{k} +
\hat{G}_{x^1} \big{]}dt + \hat{q}^{k}dW_{t}^{k},\\
\hat{p}(t,0,\xi)=0,\ \hat{p}(T,x^1,\xi)=\hat{\phi}(x^1,\xi),\ \xi \in
\mathbb{R}^{d-1},
\end{array}\right.
\end{equation}
where
$$ \tilde{b}^{1} = i\sum_{j\geq 2}a^{1j}\xi^{j},\quad
\tilde{c}=\sum_{j,k\geq 2}a^{jk}\xi^{j}\xi^{k},\quad
\tilde{\nu}^{k}=i\sum_{j\geq 2}\sigma^{jk}\xi^{j}. $$ In order to
apply Theorem \ref{thm:c1}, we define
\begin{eqnarray}\label{con:e4}
\nonumber && a=\sqrt{a^{11}},\quad \lambda^{2}=\sum_{j,k\geq
2}a^{jk}\xi^{j}\xi^{k},\quad \varepsilon =
(a\lambda)^{-1}\tilde{b}^{1},\\
&& \delta^{k}=a^{-1}\sigma^{1k},\quad
\gamma^{k}=\lambda^{-1}\tilde{\nu}^{k}.
\end{eqnarray}
From the condition \eqref{con:e1}, it is not hard to show that the condition
\eqref{con:c1} of Theorem \ref{thm:c1} is satisfied with
$\lambda_{0}=|\xi|^{2}$. However, the verification of condition
\eqref{con:c2} is rather delicate and will be put in Step 3.

Now denote $Q_{1}=\{ (t,x^1):t\in (0,T),x^1\in (0,\infty)\}$. Applying
\eqref{leq:c6} to Eq.\eqref{eq:e2}, for $0\leq r\leq m\leq n$, we have
\begin{eqnarray*}
\begin{split}
& \li (x^1)^{m}\frac{\partial^{r}}{(\partial
x^1)^{r}}(x^1\hat{p}_{x^1x^1}) \ri _{Q_1}^{2} + |\xi|^{2} \li
(x^1)^{m}\frac{\partial^{r}}{(\partial
x^1)^{r}}(x^1\hat{p}_{x^1}) \ri _{Q_1}^{2} \\
& + \li (x^1)^{m}\frac{\partial^{r}}{(\partial
x^1)^{r}}(x^1\hat{q}_{x^1}) \ri _{Q_1}^{2} + E\sup_{t\leq
T}\bigg{\|} (x^1)^{m}\frac{\partial^{r}}{(\partial
x^1)^{r}}(x^1\hat{p}_{x^1}) \bigg{\|}_{\mathbb{R}_{+}}^{2}(t) \\
\leq \ & C \bigg{\{} \li (x^1)^{m}\frac{\partial^{r}}{(\partial
x^1)^{r}}(x^1\hat{F}) \ri _{Q_1}^{2} + E\bigg{\|}
(x^1)^{m}\frac{\partial^{r}}{(\partial
x^1)^{r}}(x^1\hat{\phi}_{x^1}) \bigg{\|} _{\mathbb{R}_{+}}^{2} \bigg{\}} \\
& + C |\xi|^{-2(m-r)} \bigg{\{} \li (x^1)^{m}\frac{\partial^{m}}{(\partial
x^1)^{m}}(x^1\hat{F}) \ri _{Q_1}^{2} + E \bigg{\|}
(x^1)^{m}\frac{\partial^{m}}{(\partial x^1)^{m}}(x^1\hat{\phi}_{x^1})
\bigg{\|} _{\mathbb{R}_{+}}^{2} \bigg{\}},
\end{split}
\end{eqnarray*}
for almost all $\xi=(\xi^{2},\dots,\xi^{d})$. We multiply this inequality by
$|\xi|^{2(m-r)}$, integrate with respect to $\xi \in \mathbb{R}^{d-1}$ and
sum the results over $r=0,\dots,m$. Then we see that
\begin{eqnarray}\label{leq:e1}
\nonumber [\|x^{1}p_{x^1x^1}\|]^{2}_{m,Q}+[\|x^{1}p_{x^1y}\|]^{2}_{m,Q}
+[\|x^{1}q_{x^1}\|]^{2}_{m,Q}+E\sup_{t\leq T}[|x^{1}p_{x^1}|]^{2}_{m,\mathbb{R}^{d}_{+}}(t)\\
\leq\ C\{
[\|x^{1}F\|]^{2}_{m,Q}+[\|x^{1}\phi_{x^1}\|]^{2}_{m,\mathbb{R}^{d}_{+}} \},
\end{eqnarray}
where the subscript $y$ stands for any first-order derivative with
respect to $x^2,\dots,x^d$.

For $0\leq r\leq m\leq n$ from \eqref{leq:c7} we have
\begin{eqnarray*}
\begin{split}
& \li (x^1)^{m}\frac{\partial^{r}}{(\partial
x^1)^{r}}(x^1\hat{p}_{x^1}) \ri _{Q_1}^{2} + |\xi|^{2} \li
(x^1)^{m}\frac{\partial^{r}}{(\partial
x^1)^{r}}(x^1\hat{p}) \ri _{Q_1}^{2} \\
& + \li (x^1)^{m}\frac{\partial^{r}}{(\partial
x^1)^{r}}(x^1\hat{q}) \ri _{Q_1}^{2} + E\sup_{t\leq T}\bigg{\|}
(x^1)^{m}\frac{\partial^{r}}{(\partial
x^1)^{r}}(x^1\hat{p}) \bigg{\|}_{\mathbb{R}_{+}}^{2}(t) \\
\leq \ & C \bigg{\{} |\xi|^{-2} \li
(x^1)^{m}\frac{\partial^{r}}{(\partial x^1)^{r}}(x^1\hat{F}) \ri
_{Q_1}^{2} + E\bigg{\|} (x^1)^{m}\frac{\partial^{r}}{(\partial
x^1)^{r}}(x^1\hat{\phi}) \bigg{\|} _{\mathbb{R}_{+}}^{2} \bigg{\}} \\
& + C |\xi|^{-2(m-r+1)} \bigg{\{} \li (x^1)^{m}\frac{\partial^{m}}{(\partial
x^1)^{m}}(x^1\hat{F}) \ri _{Q_1}^{2} + E \bigg{\|}
(x^1)^{m}\frac{\partial^{m}}{(\partial x^1)^{m}}(x^1\hat{\phi}_{x^1})
\bigg{\|} _{\mathbb{R}_{+}}^{2} \bigg{\}},
\end{split}
\end{eqnarray*}
for almost all $\xi=(\xi^{2},\dots,\xi^{d})$. We multiply this inequality by
$|\xi|^{2(m-r+1)}$, integrate with respect to $\xi \in \mathbb{R}^{d-1}$ and
sum the results over $r=0,\dots,m$. Then we see that
\begin{eqnarray}\label{leq:e2}
\nonumber [\|x^{1}p_{x^1y}\|]^{2}_{m,Q}+[\|x^{1}p_{yy}\|]^{2}_{m,Q}
+[\|x^{1}q_{y}\|]^{2}_{m,Q}+E\sup_{t\leq T}[|x^{1}p_{y}|]^{2}
_{m,\mathbb{R}^{d}_{+}}(t)\\
\leq\ C\{
[\|x^{1}F\|]^{2}_{m,Q}+[\|x^{1}\phi_{x^1}\|]^{2}_{m,\mathbb{R}^{d}_{+}}
+[\|x^{1}\phi_{y}\|]^{2}_{m,\mathbb{R}^{d}_{+}} \}.
\end{eqnarray}
Combining \eqref{leq:e1} and \eqref{leq:e2}, we get that
\begin{eqnarray}\label{leq:e4}
\nonumber [\|\tilde{\psi} D^{2}p\|]^{2}_{m,Q}+E\sup_{t\leq T}[\|\tilde{\psi}
D p(t,\cdot)\|]^{2}_{m,\mathbb{R}^{d}_{+}}+[\|
\tilde{\psi} D q\|]^{2}_{m,Q}\\
\leq C([\| \tilde{\psi} F\|]^{2}_{m,Q}+[\|\tilde{\psi}
D\phi\|]^{2}_{m,\mathbb{R}^{d}_{+}}).
\end{eqnarray}

\emph{Step 2}. We now remove the assumption (A) made in Step 1.

From Lemma \ref{lem:b3} we deduce that Eq.\eqref{eq:e1} has a unique solution
$$(p,q)\in \mathbb{H}^{1}_{0}(\mathbb{R}^{d}_{+})\otimes \mathbb{H}^{0}
(\mathbb{R}^{d}_{+}),$$ such that
\begin{eqnarray*}
\interleave D p \interleave_{Q}^{2} + E\sup_{t\leq T}\| p
(t,\cdot)\|_{\mathbb{R}^{d}_{+}}^{2} + \interleave q \interleave_{Q}^{2} \leq
C (\| F \|^{2}_{\mathbb{H}^{-1}(\mathbb{R}^{d}_{+})} + \interleave \phi
\interleave_{\mathbb{R}^{d}_{+}}^{2}),
\end{eqnarray*}
where $C=C(d,\kappa,T)$. From Lemma \ref{lem:b2}, it follows that
\begin{eqnarray}\label{leq:e3}
\interleave D p \interleave_{Q}^{2} + E\sup_{t\leq T}\| p(t,\cdot)
\|_{\mathbb{R}^{d}_{+}}^{2} + \interleave q \interleave_{Q}^{2} \leq C
(\interleave \tilde{\psi} F \interleave_{Q}^{2} + \interleave \phi
\interleave_{\mathbb{R}^{d}_{+}}^{2}).
\end{eqnarray}
We shall prove that these functions are exactly what we need. Its uniqueness
is clear from the above argument. To establish the relation \eqref{leq:e},
now take an infinitely differentiable function $\zeta(y)$ defined for $y\in
\mathbb{R}$ and such that $\zeta(y)=1$ for $y\in [0,1]$, $\zeta=0$ for $y\geq
2$. Define $\zeta^{\varepsilon}(y)=\zeta(\varepsilon y)$, and
$$p^{\varepsilon}(t,x)=p(t,x)\zeta^{\varepsilon}(x^{1}),
\quad q^{\varepsilon}(t,x)=q(t,x)\zeta^{\varepsilon}(x^{1}),$$ where
the parameter $\varepsilon$ will approach infinity in the future. It
is obvious that
$$(p^{\varepsilon},q^{\varepsilon})\in \mathbb{H}^{1}_{0}(\mathbb{R}^{d}_{+})\otimes
\mathbb{H}^{0}(\mathbb{R}^{d}_{+}),$$ which satisfy the equation
\begin{equation}\label{eq:e3}
\left\{\begin{array}{l}
dp^{\varepsilon}=-[a^{ij}D_{ij}p^{\varepsilon}+\sigma^{i}D_{i}q^{\varepsilon}
+F^{\varepsilon}]dt +q^{\varepsilon}dW^k_t, \\
p^{\varepsilon}(t,x)=0,\ x\in\partial \mathbb{R}^{d}_{+},\\
p^{\varepsilon}(T,x)=\phi(x)\zeta^{\varepsilon}(x^{1}),\ x\in
\mathbb{R}^{d}_{+},
\end{array}\right.
\end{equation}
where
$$F^{\varepsilon}=F\zeta^{\varepsilon}-2\sum_{j\geq 2}
a^{1j}\zeta^{\varepsilon}_{x^1}D_{j}p^{\varepsilon} -
a^{11}\zeta^{\varepsilon}_{x^1x^1}p-\sigma^{1}
\zeta^{\varepsilon}_{x^1}q.
$$
Now let us make an assumption which will be justified later: suppose
that for an integer $r\leq n+1$ and a multi-index $\alpha$ such that
$|\alpha|\leq r$
\begin{equation}\label{con:e3}
\tilde{\psi}^{|\alpha|}D^{\alpha}p_{x},\ \tilde{\psi}^{|\alpha|}D^{\alpha}q
\in \mathbb{H}^{0}(\mathbb{R}^{d}_{+}).
\end{equation}
Note that for any
multi-index $\beta$ the functions
$$\big{|}\tilde{\psi}^{|\beta|+1}D^{\beta}\zeta^{\varepsilon}_{x^1}\big{|},\quad
\big{|}\tilde{\psi}^{|\beta|+2}D^{\beta}\zeta^{\varepsilon}_{x^1x^1}\big{|},\quad
\big{|}\tilde{\psi}^{|\beta|+1}D^{\beta}\zeta^{\varepsilon}_{x^1x^1}\big{|}$$
are bounded uniformly with respect to $\varepsilon$ and tend to zero when
$\varepsilon \rightarrow 0$. Using this and Lemma \ref{lem:b1} for $m\leq r$
and for $\varepsilon\rightarrow 0$, we get
\begin{eqnarray*}
\begin{split}
[\|\tilde{\psi} \zeta^{\varepsilon}_{x^1} p_{x}\|]_{m,Q}\ \leq \
&C\sum_{|\alpha|=m}\interleave \tilde{\psi}^{m+1}
D^{\alpha}(\zeta^{\varepsilon}_{x^1}p_{x})
\interleave_{Q}\\
\leq\ &C\sum_{|\alpha|+|\beta|=m}\interleave \tilde{\psi}^{m+1}
(D^{\beta}\zeta^{\varepsilon}_{x^1})D^{\alpha}p_{x}
\interleave_{Q}\\
=\ &C \sum_{|\alpha|+|\beta|=m}\interleave \tilde{\psi}^{|\beta|+1}
(D^{\beta}\zeta^{\varepsilon}_{x^1}) \tilde{\psi}^{|\alpha|}D^{\alpha}p_{x}
\interleave_{Q}
\rightarrow 0,\\
[\|\tilde{\psi} \zeta^{\varepsilon}_{x^1x^1} p\|]_{m,Q}\ \leq\ &C
\sum_{|\alpha|+|\beta|=m}\interleave \tilde{\psi}^{m+1}
(D^{\beta}\zeta^{\varepsilon}_{x^1x^1})D^{\alpha}p \interleave_{Q}\\
=\ &C\sum_{|\alpha|+|\beta|=m,\alpha\neq 0}\interleave
\tilde{\psi}^{|\beta|+2} (D^{\beta}\zeta^{\varepsilon}_{x^1x^1})
\tilde{\psi}^{|\alpha|-1}D^{\alpha}p
\interleave_{Q}\\
& + C \sum_{|\beta|=m}\interleave
\tilde{\psi}^{m+1}(D^{\beta}\zeta^{\varepsilon}_{x^1x^1})p \interleave_{Q}
\rightarrow 0, \\
[\|\tilde{\psi} \zeta^{\varepsilon}_{x^1} q\|]_{m,Q}\ \leq\ &C
\sum_{|\alpha|+|\beta|=m}\interleave \tilde{\psi}^{|\beta|+1}
(D^{\beta}\zeta^{\varepsilon}_{x^1}) \tilde{\psi}^{|\alpha|}D^{\alpha}q
\interleave_{Q} \rightarrow 0.
\end{split}
\end{eqnarray*}
This along with the result in Step 1 applied to Eq.\eqref{eq:e3}
after passing to the limit when $\varepsilon\rightarrow \infty$
yields
\begin{eqnarray*}
[\|\tilde{\psi} D^{2}p\|]^{2}_{m,Q}+E\sup_{t\leq T}[\|\tilde{\psi} D
p(t,\cdot)\|]^{2}_{m,\mathbb{R}^{d}_{+}}+[\|
\tilde{\psi} D q\|]^{2}_{m,Q}\\
\nonumber \leq C([\| \tilde{\psi} F\|]^{2}_{m,Q}+[\|\tilde{\psi}
D\phi\|]^{2}_{m,\mathbb{R}^{d}_{+}}),
\end{eqnarray*}
for any $m\leq r\wedge n$. This inequality implies that relation
\eqref{con:e3} with $r+1$ holds true under the assumption \eqref{con:e3} with
$r$. Since for $r=0$ relation \eqref{con:e3} follows from \eqref{leq:e3}, by
induction we obtain \eqref{leq:e}.\medskip

\emph{Step 3}. Verification of condition \eqref{con:c2}.

We need the following lemma from linear algebra.

\begin{lem}\label{lem:e1}
Suppose $K>0$, then there exists a positive number $\lambda$ such that
$$|z|^{2}+|w-A z|^{2} \geq \lambda(|z|^{2}+|w|^{2})$$
holds for any $z \in \mathbb{C}^m$, $w \in \mathbb{C}^n$ and $A \in
\mathbb{C}^{n \times m}$, so long as $|A|^2 \leq K$.
\end{lem}

\begin{proof}
Since $|Az|^{2}\leq K|z|^{2}$, it follows that for any number
$\lambda\in(0,1/2]$
\begin{eqnarray*}
\begin{split}
& |z|^{2}+|w-A z|^{2}\\
=\ &|z|^{2}+|w|^{2}- 2Re(w^{*}A z)+|A z|^{2}\\
\geq\ &|z|^{2}+|w|^{2}+|A z|^{2}-\bigg{[}(1-\lambda)|w|^{2}
+\frac{1}{1-\lambda}|A z|^{2}\bigg{]}\\
\geq\ &|z|^{2}+\lambda|w|^{2}- 2\lambda|A z|^{2}\\
\geq\ &(1-2\lambda K)|z|^{2}+\lambda|w|^{2}.
\end{split}
\end{eqnarray*}
The lemma is proved by taking $\lambda=\min\big{\{}1/2,(4K)^{-1}\big{\}}$.
\end{proof}

We now prove Theorem \ref{thm:e1}. Take any complex numbers $u,v,z_{k}
(k=1,\dots,d_{1})$ and define $\zeta=a^{-1}u,\eta=\lambda^{-1}v$. Put the
$d_{1}$-dimensional vectors
\begin{eqnarray*}
&&z=(z_{1},\dots,z_{d_1}),\\
&&\delta=(\delta^{1},\dots,\delta^{d_1}),\quad
\gamma=(\gamma^{1},\dots,\gamma^{d_1}),\\
&&\sigma^{j}=(\sigma^{j1},\dots,\sigma^{jd_1}),\ j=1,\dots,d.
\end{eqnarray*}
Recalling \eqref{con:e4} and using the standard technique in linear algebra,
we have
\begin{eqnarray*}
\begin{split}
&2|u|^{2}+2|v|^{2}+|z|^{2}-4Re(\varepsilon \bar{u}v)
-2Re(\delta \bar{u}z)-2Re(\gamma \bar{v}z)\\
=\ &2a^{2}|\zeta|^{2}+2\lambda^{2}|\eta|^{2}+|z|^{2}-4a\lambda
Re(\varepsilon \bar{\zeta}\eta)-2aRe(\delta \bar{\zeta}z)
-2\lambda Re(\gamma \bar{\eta}z)\\
=\ &
\begin{pmatrix}
\zeta & \eta & z
\end{pmatrix}
\begin{pmatrix}
2a^{11} & -2i\sum\limits_{j\geq 2}a^{1j}\xi^{j} & \sigma^{1}\\
2i\sum\limits_{j\geq 2}a^{1j}\xi^{j} & 2\sum\limits_{j,l\geq 2}a^{jl}
\xi^{j}\xi^{l} & i\sum\limits_{j\geq 2}\sigma^{j}\xi^{j}\\
(\sigma^{1})^{T} & -i\big{(}\sum\limits_{j\geq 2}\sigma^{j}\xi^{j}\big{)}^{T}
& I_{d_1}
\end{pmatrix}
\begin{pmatrix}
\bar{\zeta} \\ \bar{\eta} \\ z^{*}
\end{pmatrix}\\
=\ &\begin{pmatrix}
\zeta & \eta
\end{pmatrix}
\begin{pmatrix}
2a^{11}-|\sigma^{1}|^{2} & -i\sum\limits_{j\geq 2}(2a^{1j}
-\sum\limits_{k}\sigma^{1k}\sigma^{jk})\xi^{j}\\
i\sum\limits_{j\geq 2}(2a^{1j}-\sum\limits_{k}\sigma^{1k}\sigma^{jk})\xi^{j}
& 2\sum\limits_{j,l\geq 2}(2a^{jl}-\sum\limits_{k}\sigma^{jk}\sigma^{lk})
\xi^{j}\xi^{l}
\end{pmatrix}
\begin{pmatrix}
\bar{\zeta} \\ \bar{\eta}
\end{pmatrix}\\
& +\sum_{k=1}^{d_1}\bigg{|}z_{k}+\sigma^{1k}\zeta
+i\sum_{j\geq 2}\sigma^{jk}\xi^{j}\eta \bigg{|}^{2}\\
=\ &\sum_{j,l}b^{jl}\eta^{j}\eta^{l}+\sum_{k=1}^{d_1}\bigg{|}z_{k}
+\sigma^{1k}\zeta +i\sum_{j\geq 2}\sigma^{jk}\xi^{j}\eta \bigg{|}^{2},
\end{split}
\end{eqnarray*}
where the matrix $$ B\triangleq\big{(}b^{jl}\big{)}
=2\big{(}a^{jl}\big{)}-\sigma\sigma^{*}\geq \kappa I $$ and $$
\eta^{1}=\zeta,\eta^{j}=\xi^{j}\eta,j=2,\dots,d. $$ Note that
$\sigma\sigma^{*}\leq \kappa^{-1} I$. These along with Lemma \ref{lem:e1}
(put $A = \sigma\sigma^{*}$) yield that there exists a positive number
$\mu=\mu(\kappa)$ such that
\begin{eqnarray*}
\begin{split}
&2|u|^{2}+2|v|^{2}+|z|^{2}-4Re(\varepsilon \bar{u}v)
-2Re(\delta \bar{u}z)+2Re(\gamma \bar{v}z)\\
\geq\ &\kappa \bigg{(} |\zeta|^{2} + \sum_{j\geq 2}|\xi^{j}\eta|^{2}
+\sum_{k=1}^{d_1} \big{|}z_{k} +\sigma^{1k}\zeta
+i\sum_{j\geq 2}\sigma^{jk}\xi^{j}\eta \big{|}^{2} \bigg{)}\\
\geq\ &\mu\kappa \big{(} |\zeta|^{2}+|\xi|^{2}|\eta|^{2}+|z|^{2} \big{)}\\
\geq\ &\mu\kappa^2\big{(} |u|^{2}+|v|^{2}+|z|^{2} \big{)}.
\end{split}
\end{eqnarray*}
The proof is complete.
\end{proof}

\section{The equation in the half space}

The proof of our main theorem, Theorem \ref{thm:b2}, is based on the
following

\begin{thm}\label{thm:b2}
Let $\mathcal{D}=\mathbb{R}^{d}_{+}$ in Assumptions \ref{ass:b1} and
\ref{ass:b3}. Replace $\psi(x)$ by $\tilde{\psi}(x)=x^{1}$ in Assumption
\ref{ass:b3}. Consider the following simple form of the equation
\eqref{eq:b1}
\begin{equation}\label{eq:b2}
\left\{\begin{array}{l}
dp=-[a^{ij}D_{ij}p+\sigma^{ik}D_{i}q^k+F]dt+q^kdW^k_t, \\
p(t,x)=0,\ x\in\partial \mathbb{R}^{d}_{+},\\
p(T,x)=\phi(x),\ x\in \mathbb{R}^{d}_{+}.
\end{array}\right.
\end{equation}
Suppose that for a constant $\delta>0$ and for any $(\omega,t,x)$ we have
\begin{equation}\label{con:b1}
|a(t,x)-a_{0}(t)|\leq\delta,\quad |\sigma(t,x)-\sigma_{0}(t)|\leq\delta,
\end{equation}
where $a_{0}(t)$ and $\sigma_{0}(t)$ are some functions of $(t,\omega)$
satisfying Assumption \ref{ass:b1} and \ref{ass:b3}.

We assert that, under these assumptions, there exists a constant
$\delta(d,n,\kappa,T)>0$ such that if $\delta\leq\delta(d,n,\kappa,T)$ then
\begin{enumerate}
\item[(i)] Eq.\eqref{eq:b2} has a unique solution $(p,q)$ such that
$$ p \in \mathbb{H}^{1}_{0}(\mathbb{R}^{d}_{+}) ,\quad
\tilde{\psi} D^{2}p,q,\tilde{\psi} Dq \in
\mathbb{H}^{0}(\mathbb{R}^{d}_{+}).$$
\item[(ii)] For this solution and any
multi-index $\beta$ such that $|\beta|\leq n+1$, we have
\begin{eqnarray}\label{prp:b2}
\tilde{\psi}^{|\beta|}D^{\beta}(D p),\; \tilde{\psi}^{|\beta|}D^{\beta}q \in
\mathbb{H}^{0}(\mathbb{R}^{d}_{+}),\quad \tilde{\psi}^{|\beta|}D^{\beta}p\in
C([0,T],L^{2}(\mathbb{R}^{d}_{+})) \ (a.s.),
\end{eqnarray}
and moreover
\begin{eqnarray}\label{leq:b2}
\interleave Dp\interleave^{2}_{m+1,Q}+E\sup_{t\leq
T}\|p(t,\cdot)\|^{2}_{m+1,\mathbb{R}^{d}_{+}}+\interleave
q\interleave^{2}_{m+1,Q}\\
\nonumber \leq C(\interleave \tilde{\psi} F\interleave^{2}_{m,Q}+\interleave
\phi\interleave^{2}_{m+1,\mathbb{R}^{d}_{+}}),
\end{eqnarray}
where $m\leq n,Q=(0,T)\times \mathbb{R}^{d}_{+}$, and the constant $C$
depends only on $d,n,K,\kappa$ and $T$.
\end{enumerate}
\end{thm}

\begin{proof}
With the help of Theorem \ref{thm:e1} we can recursively define a sequence of
function pairs $(p_{r},q_{r}),r=1,2,\dots$ as solutions of the equations
\begin{eqnarray}\label{eq:f1}
\nonumber dp_{r}=-\big{[} a_{0}^{ij}D^{ij}p_{r}+\sigma_{0}^{i}D^{i}q_{r}
+(a^{ij}-a_{0}^{ij})D^{ij}p_{r-1}\\
+(\sigma^{i}-\sigma_{0}^{i})D^{i}q_{r-1} +F \big{]}dt + q_{r}dW_{t}
\end{eqnarray}
in $Q=(0,T)\times \mathbb{R}^{d}_{+}$ with the given boundary data and with
$p_{0}=0,q_{0}=0$. We denote by $C_{m}$ the right-hand side in \eqref{leq:e}.
Now let
$$I_{r}^{m}=[\|\tilde{\psi} D^{2}p_{r}\|]^{2}_{m,Q}+[\| \tilde{\psi} D q_{r}\|]^{2}_{m,Q},
\quad r\geq 1,\ m\leq n,$$
and from estimate \eqref{leq:e} we get
\begin{eqnarray*}\begin{split}
I_{r}^{m}&+E\sup_{t\leq T}[\|\tilde{\psi} D p_{r}(t,\cdot)\|]^{2}_{m,\mathbb{R}^{d}_{+}}\\
\leq\ &C_{m}+C\big{\{}[\|\tilde{\psi} (a^{ij}-a_{0}^{ij})D^{ij}p_{r-1}\|]^{2}
_{m,Q} + [\| \tilde{\psi}
(\sigma^{i}-\sigma_{0}^{i})D^{i}q_{r-1}\|]^{2}_{m,Q} \big{\}}.\end{split}
\end{eqnarray*}
With the aid of Assumption \ref{ass:b3}, \eqref{con:b1} and Lemma
\ref{lem:b1}, it is easy to check that
\begin{eqnarray*}
\begin{split}
&[\|\tilde{\psi} (a^{ij}-a_{0}^{ij})D^{ij}p_{r-1}\|]^{2}_{m,Q}\\
\leq\ &C\sum_{|\alpha|=m}\interleave \tilde{\psi}^{m+1} (a^{ij}-a_{0}^{ij})
D^{\alpha}(D^{ij}p_{r-1})\interleave^{2}_{Q}\\
&+C\sum_{|\alpha|+|\beta|=m,\beta\neq 0}\interleave \tilde{\psi}^{|\beta|}
(D^{\beta}a^{ij})\tilde{\psi}^{|\alpha|+1}D^{\alpha}(D^{ij}p_{r-1})
\interleave^{2}_{Q}\\
\leq\ &C\delta [\|\tilde{\psi} D^{2}p_{r-1}\|]^{2}_{m,Q}+
C\sum_{k<m}[\|\tilde{\psi} D^{2}p_{r-1}\|]^{2}_{k,Q}\\
\leq\ &C\delta I^{m}_{r-1}+C\sum_{k<m} I^{k}_{r-1}.
\end{split}
\end{eqnarray*}
Similarly,
\begin{equation*}
[\| \tilde{\psi} (\sigma^{i}-\sigma_{0}^{i})D^{i}q_{r-1}\|]^{2}_{m,Q} \leq
C\delta I^{m}_{r-1}+C\sum_{k<m} I^{k}_{r-1}.
\end{equation*}
Hence
\begin{equation}\label{leq:f1}
I^{m}_{r}+E\sup_{t\leq T}[\|\tilde{\psi} D
p_{r}\|]^{2}_{m,\mathbb{R}^{d}_{+}}(t) \leq C_{m}+C\delta I^{m}_{r-1} +
C\sum_{k<m}I^{k}_{r-1},
\end{equation}
where the constant $C$ depends only on $d,n,K$ and $\kappa$, and
$I_{0}^{m}=0$. It follows easily by induction from \eqref{leq:f1} that if
$\delta$ is small enough (depending on $C$), then $I_{r}^{m}\leq CC_{m}$ for
any $r\geq 1,m\leq n$ with $C=C(d,n,K,\kappa)$.

Moreover, recalling \eqref{leq:e3}, we obtain that
$$\interleave D p_{r} \interleave_{Q}^{2} + \interleave q_{r}
\interleave_{Q}^{2} + E\sup_{t\leq T}\| p_{r}(t,\cdot)
\|_{\mathbb{R}^{d}_{+}}^{2} \leq C C_{0} + C \delta I^{0}_{r-1} \leq C
C_{0}$$ with this constant $C=C(d,n,K,\kappa,T)$.

If we apply a similar argument to $p_{r}-p_{r-1}$ and $q_{r}-q_{r-1}$, we
will see that the expression
\begin{eqnarray*}
[\|D(p_{r}-p_{r-1})\|]^{2}_{m,Q}+[\| q_{r}-q_{r-1}\|] ^{2}_{m,Q}+E\sup_{t\leq
T}[\|(p_{r}-p_{r-1})(t,\cdot)\|]^{2}_{m,\mathbb{R}^{d}_{+}} \rightarrow 0
\end{eqnarray*}
as $r\rightarrow \infty$.

This obviously gives us a function pair $(p,q)$ as the limit of
$(p_{r},q_{r})$, with properties \eqref{prp:e1} and \eqref{leq:e}.
We will show this pair is what we need.

From Theorem \ref{thm:e1}, the following equation
\begin{eqnarray*}
du=-\big{[} a_{0}^{ij}D^{ij}u+\sigma_{0}^{i}D^{i}v
+(a^{ij}-a_{0}^{ij})D^{ij}p
+(\sigma^{i}-\sigma_{0}^{i})D^{i}q
+F \big{]}dt + v dW_{t}
\end{eqnarray*}
in $Q=(0,T)\times \mathbb{R}^{d}_{+}$ with the boundary data
$u(t,x)=0,x\in\partial \mathbb{R}^{d}_{+}; u(T,x)=\phi(x),x\in
\mathbb{R}^{d}_{+}$, has a unique solution $(u,v)$ with properties
\eqref{prp:e1} and \eqref{leq:e}. Applying a similar argument to $u-p_{r}$
and $v-q_{r}$, we obtain that
\begin{eqnarray*}\begin{split}
[\|\tilde{\psi} D^{2}(u-p_{r})\|]^{2}_{m,Q}+[\| \tilde{\psi} D(v-q_{r})\|]
^{2}_{m,Q}+ & E\sup_{t\leq T}[\|\tilde{\psi} D(u-p_{r})(t,\cdot)\|]^{2}
_{m,\mathbb{R}^{d}_{+}}\\
&\leq\ C\delta \interleave\tilde{\psi} D^{2}(p-p_{r-1})\interleave^{2}_{m,Q},\\
\interleave D (u-p_{r}) \interleave_{Q}^{2} + \interleave v-q_{r}
\interleave_{Q}^{2} + E\sup_{t\leq T}\| (u & - p_{r})(t,\cdot)
\|_{\mathbb{R}^{d}_{+}}^{2}\\
&\leq\ C\delta \interleave\tilde{\psi} D^{2}(p-p_{r-1})\interleave^{2}_{Q}.
\end{split}\end{eqnarray*}
By taking $r\rightarrow \infty$, it is easy to show that $u=p,v=q$, which implies
that $(p,q)$ is a solution pair of Eq.\eqref{eq:b2}.

To prove the uniqueness, assume that $(p_1,q_1),(p_2,q_2)$ are two solutions
of Eq.\eqref{eq:b2} such that
$$p_{i},Dp_{i},\tilde{\psi} D^{2}p_{i},q_{i}, \tilde{\psi} Dq_{i}\in
\mathbb{H}^{0}(\mathbb{R}^{d}_{+}).$$ It is not hard to obtain that
\begin{eqnarray*}
\interleave\tilde{\psi} D^{2}(p_1-p_2)\interleave^{2}_{Q} +\interleave D
(p_1-p_2) \interleave_{Q}^{2}
+ E\sup_{t\leq T}\| (p_1-p_2)(t,\cdot) \|_{\mathbb{R}^{d}_{+}}^{2}\\
+\interleave q_1-q_2 \interleave_{Q}^{2} \leq\ C\delta
\interleave\tilde{\psi} D^{2}(p_1-p_2)\interleave^{2}_{Q},
\end{eqnarray*}
which implies that $p_1=p_2,q_1=q_2$ by taking $\delta$ small enough
and the theorem is proved.

\end{proof}

\section{Proof of Theorem \ref{thm:b1}}

The proof is quite standard, which is divided into two steps. \medskip

\emph{Step 1}. First we will show, without using Theorem \ref{thm:b2}, that
the second inclusion in \eqref{prp:b1} holds for any solution $(p,q)$ such
that $(p,q)\in \mathbb{H}^{1}_{0}(\mathcal{D})\otimes
\mathbb{H}^{0}(\mathcal{D})$ and the first inclusion in \eqref{prp:b1} holds.
Note that for any $|\beta|\leq n+1$
\begin{eqnarray}\label{eq:g}
\begin{split}
d(\psi^{|\beta|}D^{\beta}p)=&-\psi^{|\beta|}D^{\beta}[a^{ij}D_{ij}p
+b^iD_{i}p-c p+\sigma^{ik}D_{i}q^k +\nu^kq^k+F]dt \\
&\ +\psi^{|\beta|}D^{\beta}q^kdW^k_t.
\end{split}\end{eqnarray}
By virtue of the first inclusion in \eqref{prp:b1} and
by Lemma \ref{lem:b2}(b) we have
\begin{eqnarray}\label{rel:g}
\begin{split}
& \psi^{|\beta|}D^{\beta}p \in \mathbb{H}^{1}_{0}(\mathcal{D}),\quad
\psi^{|\beta|}D^{\beta}q \in \mathbb{H}^{0}(\mathcal{D}),\\
& \psi^{|\alpha|+1}D^{\alpha}[a^{ij}D_{ij}p
+b^iD_{i}p-c p+\sigma^{ik}D_{i}q^k +\nu^kq^k+F]
\in \mathbb{H}^{0}(\mathcal{D}),
\end{split}
\end{eqnarray}
for any $|\alpha|\leq n,|\beta|\leq n+1$. In order to apply Lemma
\ref{lem:b5}, it remains to verify
\begin{equation}\label{rel:g1}
\psi^{|\beta|}D^{\beta}[a^{ij}D_{ij}p
+b^iD_{i}p-c p+\sigma^{ik}D_{i}q^k +\nu^kq^k+F]
\in \mathbb{H}^{-1}(\mathcal{D}),
\end{equation}
which holds true for $\beta =0$ as a result of Lemma \ref{lem:b2}(d) and
the observation that
$$\psi a^{ij}D_{ij}p,\ b^iD_{i}p,\ \psi c p,\ \psi\sigma^{ik}D_{i}q^k,\
\psi \nu^kq^k,\ \psi F\in \mathbb{H}^{0}(\mathcal{D}).$$
In fact, relation
\eqref{rel:g1} holds true for any $\beta$ with $|\beta|\leq n+1$. Indeed, if
$\beta=\alpha +\gamma$ with $|\alpha|\leq n$ and $|\gamma|=1$, then
\begin{eqnarray*}
\begin{split}
D^{\gamma}&(\psi^{|\alpha|+1}D^{\alpha}[a^{ij}D_{ij}p
+b^iD_{i}p-c p+\sigma^{ik}D_{i}q^k +\nu^kq^k+F])\\
=& \psi^{|\beta|}D^{\beta}[a^{ij}D_{ij}p
+b^iD_{i}p-c p+\sigma^{ik}D_{i}q^k +\nu^kq^k+F]\\
&+|\beta|(D^{\gamma}\psi)\psi^{|\alpha|}D^{\alpha}[a^{ij}D_{ij}p
+b^iD_{i}p-c p+\sigma^{ik}D_{i}q^k +\nu^kq^k+F],
\end{split}
\end{eqnarray*}
where the left-hand side is in
$\mathbb{H}^{0}(\mathcal{D})$ owing to \eqref{rel:g}.
An obvious induction gives \eqref{rel:g1} for any $|\beta|\leq n+1$.

To sum up the above arguments the second inclusion in \eqref{prp:b1} follows
from \eqref{eq:g} and Lemma \ref{lem:b5}. \medskip

\emph{Step 2.} Since $\psi F \in \mathbb{H}^{0}(\mathcal{D})$ in our
assumption, it follows from Lemma \ref{lem:b2}(d) that $F \in
\mathbb{H}^{-1}(\mathcal{D})$. Since $Da$ and $D\sigma$ are bounded,
Eq.\eqref{eq:b1} can be rewritten in divergence form (like \eqref{eq:a1}).
Therefore, by Lemma \ref{lem:b3} there exists a unique solution $(p,q)$ of
Eq.\eqref{eq:b1} belonging to $\mathbb{H}^{1}_{0}(\mathcal{D})\otimes
\mathbb{H}^{0}(\mathcal{D})$ such that
\begin{equation}\label{leq:g}
\interleave p_{x} \interleave^{2}_{Q} +\interleave q \interleave^{2}_{Q}
\leq\ C(K,L,\kappa,T) \big{(} \interleave \psi F \interleave^{2}_{Q}
+\interleave \phi \interleave^{2}_{\mathcal{D}}\big{)}.
\end{equation} \smallskip

Now take a small $\rho \in (0,\rho_{0})$ satisfying the following two
conditions.
\begin{enumerate}
\item[(1)] For any $(\omega,t)$ and $x,y\in \mathcal{D}$,
\begin{equation}\label{prp:g2}
|a(t,x)-a(t,y)|\leq \delta,\quad |\sigma(t,x)-\sigma(t,y)|\leq \delta
\end{equation}
if $|x-y|\leq \rho$, where $\delta=\delta(d,n,\kappa,T)$ is taken from
Theorem \ref{thm:b2}.

Here we only require the equicontinuity of $a(\omega,t,\cdot)$ and
$\sigma(\omega,t,\cdot)$ on $\bar{\mathcal{D}}$ for all
$(\omega,t)\in\Omega\times[0,T]$, which is implied by Assumption
\ref{ass:b4}.

\item[(2)] If $x,y$ belong to the same domain $U$ from Assumption
\ref{ass:b2}, then for any $(\omega,t)$,
\begin{eqnarray}\label{prp:g1}
\nonumber &\big{|} D\Psi(x) a(t,x) (D\Psi(x))^{*}- D\Psi(y) a(t,y)
(D\Psi(y))^{*}\big{|}\leq \delta_1, &\\
&\big{|} D\Psi(x) \delta(t,x) - D\Psi(y) \delta(t,y) \big{|} \leq \delta_1&
\end{eqnarray}
if $|x-y|\leq \rho$, where the constant $\delta_1 =
\delta_{1}(d,n,\kappa^2,T)$ and $\Psi$ is the inverse for $\Phi$.
\end{enumerate}

Take any functions $\zeta,\eta \in C_{0}^{\infty}(\mathbb{R}^{d})$ such that
$\textrm{supp}(\zeta)\subset B_{2\rho}(0)$, $\int \zeta dx = 1$, $\zeta(x)=1$
for $|x|\leq \rho$, and $\eta(y)=1$ for $|y|\leq 1$, $\eta(y)=0$ for $|y|\geq
2$, $0\leq \eta \leq 1$. For any $z\in \mathbb{R}^d$ define
$$\zeta^{z}(x)=\zeta(x-z),\quad p^{z}(t,x)=p(t,x)\zeta^{z}(x),\quad
q^{z}(t,x)=q(t,x)\zeta^{z}(x).$$ Now consider two case. \medskip
\\
\emph{Case 1.} $\textrm{dist}(z,\partial \mathcal{D})\leq 2\rho_0$. Then
$\mathcal{D}\cap \textrm{supp}(p^{z})\subset U \cap \mathcal{D}$, where $U$
is taken from Assumption \ref{ass:b2}. Define
$$u^{z}(t,y)=p^{z}(t,x)\eta(y),\quad v^{z}(t,y)=q^{z}(t,x)\eta(y),$$
with $x=\Phi(y)$ in $\tilde{Q}=(0,T)\times \mathbb{R}^{d}_{y,+}$. Obviously
$(u^{z},z^{z})\in \mathbb{H}^{1}_{0}(\mathbb{R}^d_{+})\otimes
\mathbb{H}^{0}(\mathbb{R}^d_{+})$. It is not hard to check that the functions
$u^{z},z^{z}$ satisfy the equation
\begin{equation}\label{eq:g1}
du^{z} = -[\tilde{a}^{ij}u^{z}_{y^iy^j}+\tilde{\sigma}^{ik} v^{z,k}_{y^i} +
\tilde{F}]dt + v^{z,k}dW_{t}^{k},
\end{equation}
where (observe that $u^{z}=0,v^{z}=0$ whenever $\eta \neq 1$)
\begin{eqnarray*}\begin{split}
& x=\Phi(y),\ x_0=\Phi(0),\ L^{0}=a^{ij}\partial^{2}_{x^ix^j}
+b^{i}\partial_{x^i},\\
& \tilde{a}^{ij}(t,y) = a^{rs}(t,x)\Psi^{i}_{x^r}\Psi^{j}_{x^s}(x)\eta(y)
+ a^{rs}(t,x_0)\Psi^{i}_{x^r}\Psi^{j}_{x^s}(x_0)(1-\eta(y)),\\
& \tilde{\sigma}^{ik}(t,y) = \sigma^{rk}(t,x)\Psi^{i}_{x^r}(x)\eta(y)
+\sigma^{rk}(t,x_0)\Psi^{i}_{x^r}(x_0)(1-\eta(y)),\\
& \tilde{F}(t,y) = (\zeta^{z}F)(t,x)\eta(y)+p_{x^r}(t,x)\Theta_{1}^{r}(t,y)
+p(t,x)\Theta_{2}(t,y)+q^{k}(t,x)\Theta_{3}^{k}(t,y),\\
& \Theta_{1}^{r}(t,y) = [(\zeta^{z}L^{0}\Psi^{i})(t,x)\Phi^{r}_{y^i}(t,y)
-2(a^{ij}\zeta^{z})(t,x)]\eta(y),\\
& \Theta_{2}(t,y) = [(\zeta^{z}_{x^r}L^{0}\Psi^{i})(t,x)\Phi^{r}_{y^i}(t,y)
-(L^{0}\zeta^{z}+c\zeta^{z})(t,x)]\eta(y),\\
& \Theta_{3}^{k}(t,y) = (\nu^{k}\zeta^{z}-\sigma^{rk}\zeta^{z}_{x^r})
(t,x)\eta(y).
\end{split}\end{eqnarray*}
The choice of $\rho$ (see \eqref{prp:g1}) and our assumptions about
$|(D\Phi)\zeta|$ from Assumption \ref{ass:b2} show that
$\tilde{a},\tilde{\sigma}$ satisfy the first condition in Assumption
\ref{ass:b3} and condition \eqref{con:b1} with $\kappa^{2}$ instead of
$\kappa$, with the corresponding $\delta$, and with $a_0(t)=\tilde{a}(t,0)$
(see Theorem \ref{thm:b2}). It is also easy to check that other conditions of
Assumption \ref{ass:b3} with $n=0$ (and another constant $K$) are satisfied
for Eq.\eqref{eq:g1}. Therefore, by Theorem \ref{thm:b2} with $n=0$ for
$\tilde{F}$ as above, Eq.\eqref{eq:g1} has a unique solution $(u,v)$ such
that
$$ u \in \mathbb{H}^{1}_{0}(\mathbb{R}^{d}_{+}) ,\quad \tilde{\psi} u_{yy},v,
\tilde{\psi} v_{y} \in \mathbb{H}^{0}(\mathbb{R}^{d}_{+}),$$ with
$\tilde{\psi}(y)=y^1$. Moreover, this solution is unique in the space
$\mathbb{H}^{1}_{0}(\mathbb{R}^{d}_{+})\otimes
\mathbb{H}^{0}(\mathbb{R}^{d}_{+})$, since Eq.\eqref{eq:g1} can be rewritten
in divergence form. Hence we deduce that
$$u^{z},u_{y}^{z},\tilde{\psi}u^{z}_{yy},v^{z},\tilde{\psi} v^{z}_{y}
\in \mathbb{H}^{0}(\mathbb{R}^{d}_{+}).$$

Now we apply the estimates from Theorem \ref{thm:b2} to Eq.\eqref{eq:g1}.
First of all, we note that for any function $h$
\begin{equation}\label{for:g}
D^{\alpha}_{y}(\partial_{y^i}h) =D^{\alpha}_{y}(h_{x^{r}}\Phi^{r}_{y^i}) =
\sum_{\beta+\gamma =\alpha} c^{\alpha}_{\beta\gamma}(D^{\beta}_{y}h_{x^{r}})
D^{\gamma}_{y}\Phi^{r}_{y^i},
\end{equation}
where $c^{\alpha}_{\beta\gamma}$ are some constants. Observe that (recalling
$\tilde{\psi}=y^1$) from our assumptions on $\psi$ and $\Phi$ there exists a
constant $C$ depending only on $\kappa$ such that
$$C^{-1}\tilde{\psi}(y)\leq \psi(x)\leq C\tilde{\psi}(y)$$
with $x=\Phi(y),y\in 2B_{+}$. Hence from the formula \eqref{for:g} and
Assumption \ref{ass:b3} , we get
$$\tilde{\psi}^{|\alpha|} \big{(} |D^{\alpha}_{y}\Theta_{1}|
+|D^{\alpha}_{y}(\tilde{\psi} \Theta_{2})|+ |D^{\alpha}_{y}(\tilde{\psi}
\Theta_{3})| \big{)} \leq C, \quad |\alpha|\leq n.$$ Denote
\begin{eqnarray*}\begin{split}
&\tilde{Q}=(0,T)\times \mathbb{R}^{d}_{y,+}, &D(z,r)=
B_{r}(z)\cap\mathcal{D},\\ &Q(z,r)=(0,T)\times D(z,r),~~~
&\rho_{1}=8K\rho_{0},~~~ m\leq n.\end{split}
\end{eqnarray*}
From Lemma \ref{lem:b1}(b) we have (functions $\eta=\eta(y)$ are different in
different place)
\begin{eqnarray*}
\begin{split}
[\|\tilde{\psi}\zeta^{z}F\eta\|]^{2}_{m,\tilde{Q}}\ =\ & [\|\psi
F(\tilde{\psi}\psi^{-1}\zeta^{z}\eta)\|]^{2}_{m,\tilde{Q}} \leq
C\sum_{|\alpha|\leq m}\interleave \eta^{\alpha}\tilde{\psi}^{m}
D^{\alpha}_{y}(\psi F) \interleave ^{2}_{\tilde{Q}}\\
\leq\ &C \sum_{|\alpha|\leq m}\interleave \eta^{\alpha}\psi^{m}
D^{\alpha}(\psi F) \interleave ^{2}_{Q} \leq C \sum_{|\alpha|\leq
m}\interleave \psi^{|\alpha|}
D^{\alpha}(\psi F) \interleave ^{2}_{Q(z,\rho_1)}\\
=\ & C \interleave \psi F \interleave^{2}_{m,Q(z,\rho_1)},
\end{split}
\end{eqnarray*}
where $\eta^{\alpha}=\eta^{\alpha}(y)$ are some bounded functions vanishing
if $|y|\geq 2$. Moreover, using a formula similar to \eqref{for:g} with $x$
and $y$ interchanged and Lemma \ref{lem:b1} (a),(b), we have
\begin{eqnarray*}
\begin{split}
[\|p_x\|]^{2}_{m+1,Q(z,\rho)}\ \leq\ & [\|(\zeta^{z}p)_x\|]^{2}_{m+1,Q} = C
\sum_{|\alpha|=m+1}
\interleave \psi^{m+1}D^{\alpha}(\zeta^{z}p\eta) \interleave ^{2}_{Q}\\
=\ & C \sum_{|\alpha|=m+1}\interleave \psi^{m+1}D^{\alpha}u^{z}_{x}
\interleave ^{2}_{Q} \leq C \sum_{|\alpha|\leq m+1}\interleave
\tilde{\psi}^{m+1}
D^{\alpha}_{y}u^{z}_{y} \interleave ^{2}_{\tilde{Q}}\\
\leq \ &C \sum_{r\leq m} [\|\tilde{\psi}u^{z}_{yy}\|]^{2}_{r,\tilde{Q}}
= C \interleave \tilde{\psi}u^{z}_{yy} \interleave^{2}_{m,\tilde{Q}},\\
[|p(t,\cdot)|]^{2}_{m+1,D(z,\rho)}\ \leq\ &
[|\zeta^{z}p(t,\cdot)|]^{2}_{m+1,\mathcal{D}} \ \leq\
C \|\tilde{\psi}u^{z}_{y}(t,\cdot)\|^{2}_{m,\mathbb{R}^{d}_{y,+}},\\
[\|q\|]^{2}_{m+1,Q(z,\rho)}\ \leq\ & [\|\zeta^{z}q\|]^{2}_{m+1,Q} \ \leq\ C
\interleave \tilde{\psi}v^{z}_{y} \interleave^{2}_{m,\tilde{Q}}.
\end{split}
\end{eqnarray*}
Bearing in mind several similar estimates and from Theorem \ref{thm:b2} we
obtain
\begin{eqnarray}\label{leq:g1}
\begin{split}
&\interleave p_{x} \interleave^{2}_{m+1,Q(z,\rho)} +E\sup_{t\leq
T}\|p(t,\cdot)\|^{2}_{m+1,D(z,\rho)}
+\interleave q \interleave^{2}_{m+1,Q(z,\rho)}\\
\leq\ & C \big{(} \interleave \psi F \interleave^{2}_{m,Q(z,\rho_1)}
+\interleave \phi \interleave^{2}_{m+1,D(z,\rho_1)} +\interleave p
\interleave^{2}_{m+1,Q(z,\rho_1)} +  \interleave q
\interleave^{2}_{m,Q(z,\rho_1)}\big{)}
\end{split}
\end{eqnarray}
for any $m\leq n$ under the same condition of finiteness. \medskip
\\
\emph{Case 2.} $\textrm{dist}(z,\partial \mathcal{D})\geq 2\rho_0$. This case
can easily be reduced to the first one. Indeed, we can replace the domain
$\mathcal{D}$ by any half space with boundary lying at a distance $2\rho_0$
from $z$. In this situation it is unnecessary to flatten the boundary and to
make any change of coordinates. Then as above we get an estimate similar to
\eqref{leq:g1} with norms defined with the help of the distance from the new
boundary. As above from this estimate we get \eqref{leq:g1}, keeping in mind
that the new distance and $\psi$ are bounded away from zero on
$\textrm{supp}(\zeta^{z}p)$.

Integrating \eqref{leq:g1} with respect to all $z\in \mathbb{R}^{d}$, we
obtain that
\begin{eqnarray}\label{leq:g2}
\begin{split}
\interleave p_{x} &\interleave^{2}_{m+1,Q} +E\sup_{t\leq
T}\|p(t,\cdot)\|^{2}_{m+1,\mathcal{D}}
+\interleave q \interleave^{2}_{m+1,Q}\\
\leq\ & C \big{(} \interleave \psi F \interleave^{2}_{m,Q} +\interleave \phi
\interleave^{2}_{m+1,\mathcal{D}} +\interleave p \interleave^{2}_{m+1,Q}
+\interleave q \interleave^{2}_{m,Q}\big{)},
\end{split}
\end{eqnarray}
where the constant $C$ depends on $n,K,\rho_{0},\kappa$ and modulus of
continuity of $a$ and $\sigma$ (see \eqref{prp:g2}).

For $m=0$ the right-hand of \eqref{leq:g2} is finite since $ p,p_{x},q\in
\mathbb{H}^{0}(\mathcal{D})$. It follows from induction that the right-hand
of \eqref{leq:g2} is finite for any $m\leq n$, which indicates the first
inclusion in \eqref{prp:b1}. Observe that the above estimate also holds if we
replace the initial time zero by any $s\in[0,T)$. Then by the Gronwall
inequality and induction, we get
\begin{eqnarray}\label{leq:g3}
\begin{split}
\interleave p_{x} &\interleave^{2}_{n+1,Q} +E\sup_{t\leq
T}\|p(t,\cdot)\|^{2}_{n+1,\mathcal{D}}
+\interleave q \interleave^{2}_{n+1,Q}\\
\leq\ & C \big{(} \interleave \psi F \interleave^{2}_{n,Q} +\interleave \phi
\interleave^{2}_{n+1,\mathcal{D}} +\interleave q \interleave^{2}_{n,Q}\big)\\
\leq\ & C \big{(} \interleave \psi F \interleave^{2}_{n,Q} +\interleave \phi
\interleave^{2}_{n+1,\mathcal{D}} +\interleave q \interleave^{2}_{Q}\big{)}.
\end{split}
\end{eqnarray}
This along with the inequality \eqref{leq:g} yields the estimate
\eqref{leq:b1}.

The uniqueness of the solution follows from Lemma \ref{lem:b3}. The proof is
complete.

\begin{rmk}
Instead of Assumption \ref{ass:b4}, we could also obtain the inequality
\eqref{leq:g2} by a weaker condition, i.e., the equicontinuity of
$a(\omega,t,\cdot)$ and $\sigma(\omega,t,\cdot)$ on $\bar{\mathcal{D}}$
(recall \eqref{prp:g2}), provided the last two terms of \eqref{leq:g2} is
finite for $m=0$. However, this condition, which is natural for SPDEs (see
\cite{Kryl94}), is not enough for us to estimate the quantity $\interleave q
\interleave_{Q}$ (see \eqref{leq:g3}) which does not appear in a SPDE. It is
interesting to seek a condition, that is weaker than Assumption \ref{ass:b4},
also guarantees the finiteness of $\interleave q \interleave_{Q}$.
\end{rmk}

\bibliographystyle{plain}


\end{document}